\documentclass[12pt,leqno]{article}

\usepackage[english]{babel}
\usepackage{amssymb,amsmath}
\usepackage{stmaryrd}
\usepackage[utf8]{inputenc}
\usepackage{tikz}
\usepackage{graphicx}
\usepackage[T1]{fontenc}
\usepackage{latexsym}
\usepackage{url}
\usepackage{datetime}
\usepackage{here}
\usepackage{xcolor}

\usepackage{amsthm}
\usepackage{shadethm}	

\usepackage{hyperref}

\definecolor{shadethmcolor}{gray}{0.92}     
\newtheoremstyle{theorem}		    % name
	{}				    % Space above
	{}				    % Space below
	{}				    % Body font
	{}				    % Indent amount
	{\bfseries}			    % Thm head font
	{}				    % Punctuation after thm head
	{\newline}			    % Space after thm head: " " = normal interword space;
	{\thmname{#1}\thmnumber{~#2}\thmnote{\itshape~~#3}}	
\theoremstyle{theorem}
\newshadetheorem {theorem}                           {Theorem}  [section]
\newshadetheorem {axiom}                  [theorem]  {Axiom}
\newshadetheorem {assertion}              [theorem]  {Assertion}
\newshadetheorem {corollary}              [theorem]  {Corollary}
\newshadetheorem {definition}             [theorem]  {Definition}
\newtheorem {example}                     [theorem]  {Example}

\newshadetheorem {fact}                   [theorem]  {Fact}
\newshadetheorem {lemma}                  [theorem]  {Lemma}
\newshadetheorem {observation}                 [theorem]  {Observation}
\newshadetheorem {notation}               [theorem]  {Notation}   
\newshadetheorem {proposition}            [theorem]  {Proposition}
\newshadetheorem {assumption}             [theorem]  {Assumption}
\newtheorem {remark}                      [theorem]  {Remark}

\numberwithin{equation}{section}

\newcommand{\ignore}[1]{}

\newcommand{\nolig}{\kern0.3pt}

\newcommand{\smallbox}[1]{\parbox{\textwidth}{\small #1}}

\DeclareRobustCommand{\nec}{%
\mathbin{\text{\scalebox{.9}{\raisebox{0.0pt}{$\Box$}}}}%
}
\DeclareRobustCommand{\poss}{%
  \mathbin{\text{\scalebox{.75}{\raisebox{1.5pt}{\rotatebox[origin=c]{45}{$\Box$}}\,}}}%
}
\newcommand{\nesi}{\stackrel{{\scriptstyle\nec}}{\sim}}

\renewcommand{\qed}{\hspace*{\fill}{{\tiny$\blacksquare$}}}
\newcommand{\beginproof}{\medskip\noindent{\textit{Proof}:~}}

\newcommand{\Sp}{(A)}
\newcommand{\SpL}{(A$_L$)}

\newcommand{\Spiter}{(A$^{\!\text{iter}}$)}

\newcommand{\SpsrL}{(A$_L^{\text{self}}$)}
\newcommand{\SpsreL}{(A$_L^{\text{self},\exists}$)}
\newcommand{\SpastL}{(A$^{{\textstyle\ast}}_L$)}
\newcommand{\Spast}{(A$^{{\textstyle\ast}}$)}

\newcommand{\cT}{c_{\text{T}}}
\newcommand{\rT}{r_{\text{T}}}

\newcommand{\Ac}{\mathcal{A}}

\newcommand{\Mc}{\mathcal{M}}

\newcommand{\Oc}{\mathcal{O}}

\newcommand{\Vc}{\mathcal{V}}

\newcommand{\Mo}{{\textsf{Mo}}}
\newcommand{\Tu}{{\textsf{Tu}}}
\newcommand{\We}{{\textsf{We}}}
\newcommand{\Th}{{\textsf{Th}}}
\newcommand{\Fr}{{\textsf{Fr}}}
\newcommand{\none}{{\textsf{none}}}

\newcommand{\svdash}{\vdash_{\text S4}}
\newcommand{\sfdash}{\vdash_{\text S5}}

\newcommand{\McG}{\Mc_{\text{G}}}
\newcommand{\WG}{W_{\text{G}}}
\newcommand{\EG}{E_{\text{G}}}
\newcommand{\VG}{V_{\text{G}}}

\newcommand{\Vn}{\mathcal{V}_{{\scriptstyle\nec}}}

\begin{document}

\ignore{\title{{\Large\bf DRAFT, do not distribute}\\[24pt]
Dismantling the Surprise Test ``Paradox''\\[6pt]
(Journal version)}}

\title{Dismantling the Surprise Test ``Paradox''}
%
%Rational Choice makes
%the Surprise Test ``Paradox'' vanish

\author{Martin Dietzfelbinger\\Technische Universität Ilmenau}
\date{\today}

\maketitle

\begin{center}
\emph{In memoriam Martin Mundhenk (1961--2024).}
\end{center}

\begin{abstract}
\noindent
This paper is about the following well-known story: A teacher announces to her students 
a test for the following week, such that the day of the test will be ``surprising''. The students use
this announcement as the basis for a ``logical derivation'' leading to a contradiction, 
which they (falsely) interpret as saying that there cannot be a test.
Nonetheless, the teacher gives a test e.g. on Wednesday, (intuitively) surprising the students. --
Its curious turns give the story the flavor of a paradox.  
The story runs under several alternative names, 
e.g. the \emph{unexpected hanging paradox} and the \emph{prediction paradox}.
Discussions and analyses of the story in the philosophical and mathematical 
literature are abundant, spanning the period from the late 1940s
until today. 
Apparently, none of the many explanations that have appeared
has been generally accepted as conclusive and convincing.
This paper offers a fresh view, on the basis of propositional logic.
``Surprise'' relative to some axiom system is captured as unprovability of a certain propositional formula.
``Knowledge'' in general corresponds to axiom systems, and it can be gained by mathematical proofs.  
In this framework, the notorious property of self-reference in the announcement of the teacher can be cleanly accommodated.
All errors in the argumentation of the students are identified.
A general analysis shows that the announcement can be realized in a way that matches intuition,
and it shows how to resolve doubts regarding this intuitive view. 
This is the first mathematically precise analysis of the story of the surprise test 
showing that self-reference, full power of mathematical proofs, 
and truthfulness of the teacher can coexist without creating a contradiction.
The ``paradox'' vanishes into thin air.
In order to facilitate comparisons with the many treatments in the literature that used modal logic,
a version based on system S5 is also given.
A formula $\sigma$ is identified that formalizes ``there will be a surprising test'', 
and it is shown that the students take the announcement to mean $\nec\sigma$
while in fact the information conveyed by it is not stronger than $\poss\!\sigma$.
This dissolves all contradictions or ``paradoxical'' issues that have been seen in the announcement. 
\end{abstract}

\tableofcontents

%%%%%%%%%%%%%%%%%%%%%%%%%%%%%%%%%%%%%%%%%%%%%%%%%%%%%%%%%%%%%%%%%%%%
\section{Introduction}\label{sec:introduction}
%%%%%%%%%%%%%%%%%%%%%%%%%%%%%%%%%%%%%%%%%%%%%%%%%%%%%%%%%%%%%%%%%%%%

\subsection{The story: A paradox?}\label{sec:story:paradox}
%%%%%%%%%%%%%%%%%%%%%%%%%%%%%%%%%%%%%%%%%%%%%%%%%%%%%%%%%%%%%%%%%%%%

Some ``narrator'' tells the following story.

\medskip

\begin{small}
\begin{itemize}
	\item[\quad] One Friday, a teacher makes the following announcement in her math class. 
\begin{itemize}
\item[{\Sp}] ``On some day of next week there will be a test. And it will be surprising, that means that
on the morning of the day of the test you won't know that it is on that day.'' 
\end{itemize}
\end{itemize}
\end{small}
The narrator adds
that there is a math lesson on each day of the week
and that the rules of the school prohibit having more than one test in a week.
\begin{small}
\begin{itemize}
\item[\quad] 
The students trust their teacher, and so they are sure that she speaks the truth. 
But they are also clever, and they start thinking and deducing, as follows. 
``We know there is a test next week, and that it is surprising. 
But then it cannot be on Friday. Otherwise we would have seen
on Friday morning that the previous days have passed without test, 
we would know that the test must be on this day, and we would not be surprised.
So the test definitely is on one of the days from Monday through Thursday. 
But hey! now we can use the same argument to conclude that the test cannot be on Thursday either. 
And we can repeat this for Wednesday, Tuesday, Monday, which gives us that
the test cannot take place on any of the 
days of next week. This means there is no test, hooray!'' 
\end{itemize}
\end{small}
In their excitement, the students hardly notice that in essence
they have found out that announcement {\Sp} \emph{cannot be realized}.
\begin{small}
\begin{itemize}
\item[\quad]
When the teacher comes to class on Wednesday and hands out the test papers, the students are truly surprised
that the test is on this day.
\end{itemize}
\end{small}
The students did not foresee that the test will be on Wednesday. 
This looks as if the teacher \emph{has realized} announcement {\Sp} after all.

\bigskip

In this section, we will sketch our way of viewing the story. For remarks on the history and the extensive literature on
the subject see Section~\ref{sec:background}. 

In principle, the teacher has a choice between six options,
namely having the test on one of the five days or having no test at all.
Clearly it is irrelevant at which point in time the teacher makes her decision, 
and we can consider her choice $\rT$ as fixed, but unknown to the students.
Announcement {\Sp} becomes: ``$\rT$ is a day $d$, and the test on day $d$ is surprising.''
In this form it makes sense to ask whether {\Sp} is true or not, 
but it is important to keep in mind that the answer depends on the value $\rT$.  

The simplest form of a logical paradox is a situation in which
one can seemingly prove a statement S and can also seemingly prove the negation of S.
The first part of the story with the argument of the students looks like a proof of
``If {\Sp} is announced and is true with $\rT$, 
then there cannot be a test, let alone a surprising one.''
The last part of the story seems to say that notwithstanding the argument 
of the students the teacher can have her test and {\Sp} can be made true,
even though {\Sp} has been announced.
This situation gives the impression of a logical paradox as just characterized. 

It is easy to see that it is sufficient to concentrate on the first part. 
Namely, already a superficial check reveals a blatant logical error in the argument of the students.%
\footnote{This has been observed many times before, of course.}
They start from the assertion that choice $\rT$ will result in a surprising test 
and reach the conclusion that there cannot be a test at all.
This is not a result, but a contradiction, which the students do not recognize as such. 
Instead, they consider ``there is no test'' a proven fact.
But intermediate results in a chain of inferences that leads to a contradiction
cannot be taken as facts, since from a contradictory situation one can prove everything,
so that such a proof is meaningless.
This terrible mistake on the side of the students in combination with the commonsensical fact that
the teacher can realize {\Sp} easily raises the suspicion that something is wrong with the ``logical argument'' of the students. 
But what exactly is wrong with it? This is the version of the ``surprise test paradox (STP)'' that we will study here.

A paragraph in T. Y.~Chow's paper~\cite{Chow:1998} nicely describes the approach of the present paper. Chow writes:
\begin{quote}\small
``In general, there are two steps involved in resolving a paradox. First, one establishes
precisely \emph{what the paradoxical argument is}. Any unclear terms are defined carefully and
all assumptions and logical steps are stated clearly and explicitly, possibly in a formal
language of some kind. Second, \emph{one finds the fault in the argument}. Sometimes, simply
performing step one reveals the flaw, e.g., when the paradox hinges on confusing two
different meanings of the same word, so that pointing out the ambiguity suffices to dispel
the confusion. In other cases, however, something more needs to be done; one must locate
the bad assumptions, the bad reasoning, or (in desperate circumstances) the flaw in the
structure of logic itself.''
\end{quote}
We find a (seemingly new) formalization of STP in terms of propositional logic, 
with axiom systems, proofs, and models (which are the assignments to the propositional variables).
Simple statements like ``there is a test on Tuesday or Wednesday or Thursday'' are formulated
as plain propositional formulas. 
``Knowledge'' is represented by axiom systems, 
and ``surprise'' relative to some ``knowledge'' by the unprovability of certain formulas from the corresponding axiom system.
Gaining ``knowledge'' means establishing stronger axiom systems by mathematical argumentation. 
This setup is already sufficient to prove that the argumentation of the students is not 
mathematically sound, and to clarify which single steps are not justified, and why. 

This first analysis pinpoints the spots in the STP story 
where sloppy usage of everyday language camouflages illegal steps in the argumentation of the students.
Concretely, these are the following:  
The sentence ``\emph{There must be a test}'' (as a consequence of {\Sp} being true with $\rT$) could mean (i) either that the concrete choice $\rT$
must be one of the days
(ii) or that it is absolutely impossible that the teacher chooses to have no test.
An analogous distinction can be made for ``\emph{there cannot be a test on Friday}'' as a conclusion 
from ``it is impossible that there is no test'' and {\Sp} being true with $\rT$. 
In the mathematical formulation it will become clear that (i) and (ii) are different things;
in everyday language this distinction is easily glossed over.  
In a sense, Chow's recipe cited above is carried out literally:
We realize that ``the paradox hinges on confusing two
different meanings of the same word'', where in our case it is not a word, but a sentence.

In a second part, we will analyze the situation laid out in STP in a more general way.
A rather superficial change is that we admit ``laws'', which can definitely exclude certain options.
(With this change, all variants of STP to be found in the literature can be covered.)
A notorious feature of STP is that when {\Sp} is taken at its full strength, 
it is ``self-referential'', in saying that the students will be ``surprised'' by the test 
\emph{even if they take into account that \emph{{\Sp}} with $\rT$ is true}.
This feature was recognized in~\cite{Shaw:1958} and was made responsible for the contradiction appearing in the story in 
that paper and several subsequent ones, like~\cite{Fitch:1964,Kaplan:Montague:1960,Kritchman:Raz:2010}.
This feature of ``self-referentiality'' woven into STP can smoothly be accommodated in our framework.

Not so much as a technical ingredient as for gaining intuition 
we will point out an analogy in simple (one-move) games and use the concept of a 
\emph{rational teacher}, or a teacher who \emph{behaves rationally}. 
This is a teacher who when given a choice between several alternatives, 
some advantageous for her, some disadvantageous, will always choose an advantageous alternative, if possible.%
\footnote{An in-depth study of the interplay of games and (knowledge about) rational behavior is given in~\cite{Gintis:2014}, in the framework of epistemic logic.
In comparison, our usage of the concept is very shallow. }
We will use this concept in two ways. 

Firstly, we will see that the main ingredient in the argument of the students is
the idea that rational behavior of their teacher can be used as a mathematical truth.%
\footnote{See~\cite{Sober:1998} for the opposite approach to that explicitly make ``rational behavior of the teacher''
a postulate that can be ``known'' and used in proofs by the students.}
The combination of {\Sp} being true and the teacher behaving rationally as a mathematical fact
is contradictory. The rather obvious way out of this situation is
to recognize that if one wants to have ``surprise'' one cannot postulate at the same time
that the behavior of the teacher must be fully predictable.
Moreover, it seems to be a dubious decision to set up axioms that model a situation where the teacher's choices 
are restricted when in reality she obviously can do what she wants.
 
Secondly, we note that it can never be excluded that the teacher does act rationally, so that all 
arguments the students propose must work in particular given a rational teacher. 
In this case the statement ``My choice $\rT$ realizes {\Sp}''
does not provide more information than ``There is some $\rT$ that realizes {\Sp}''.
Analyzing this ``existential reading'' of {\Sp} shows that even taking into account the full ``self-reference'' property of {\Sp}
the announcement is much weaker than it appears in the argumentation of the students.%
\footnote{Vinogradova~\cite{2023:Vinogradova} stresses the difference between certainty and possibility, but in a technical framework totally different 
from ours, namely constructive mathematics.}
When viewed in this way, the self-reference in STP has no effect at all,
and the students do not learn anything from {\Sp} 
beyond that it is the goal of the teacher to have a surprising test.

We arrive at the first mathematically precise analysis of the story of the surprise test 
that shows that full self-reference, full power of mathematical proofs, and truthfulness of the teacher can coexist 
without creating a contradiction. 
The ``paradox'' in the surprise test ``paradox'' vanishes into thin air.

%%%%%%%%%%%%%%%%%%%%%%%%%%%%%%%%%%%%%%%%%%%%%%%%%%%%%%%%%%%%%%%%%%%
\section{Background}\label{sec:background}
%%%%%%%%%%%%%%%%%%%%%%%%%%%%%%%%%%%%%%%%%%%%%%%%%%%%%%%%%%%%%%%%%%%
The study of our story has a history of more than 80 years. 
It seems to have surfaced in the middle of the 1940s. 
For a detailed description of the origins and the literature up to 1988 see the corresponding chapter in Sorensen's book~\cite{sorensen1988blindspots}.
A pivotal point of reference is Chow's paper~\cite{Chow:1998}, which 
gives an account of the situation in 1998 from the mathematical viewpoint.
Its version updated in 2011~\cite{chow2011surprise} 
provides an extended reference list of about 200 articles and books dealing with the ``paradox'' or giving it some treatment.
More papers keep appearing since that time, among them the knowledgeable survey by Earman~\cite{Earman:2021}.

In view of the vast literature it is not possible to discuss the whole picture here. 
We just point to the origins and the main directions.
The ``paradox'' was introduced into the philosophical discussion by O'Connor~\cite{OConnor:1948} and 
Scriven~\cite{Scriven:1951}. The latter pointed out the discrepancy between the result of the ``logical derivation''
and ``reality'' and distinguishes ``announcements as ordainments'' and ``announcements as statements'',
which might be viewed as a distinction between a situation where the word of the teacher in {\Sp} 
forces her to behave rationally and the absence of such an effect. 
W.~V.~Quine~\cite{Quine:1953} denied the existence of a contradiction
and argued that the ``initial step'' of the students, in which they state that it is certain 
that there must be a test, has no justification. 
The problem with his argumentation is that the concept of ``knowing'' something 
is used without giving it a precise meaning. 
Quine also describes the variant of STP, called the ``birthday present `paradox'\kern1pt'',
where there is only one day, Friday, say. Shaw~\cite{Shaw:1958} observed the ``self-reference'' feature
in the announcement and made it responsible for the contradiction arising.
Kaplan and Montague~\cite{Kaplan:Montague:1960} go in a similar direction, with more formalism,
including operators that express ``knowledge'' at certain points in time, 
where ``knowledge'' means provability, and a mechanism that allows to express mathematical provability within formulas,
in essence based on Goedelization.%
\footnote{They introduce the ``knower paradox'', which says that the sentence ``$X$ holds, but I do not know it'' cannot be known to the speaker. 
There is a resemblance to our ``there will be a test, and you will not be able to prove when before it happens'' and the fact that
if this is known to be a fact to the students, contradiction ensues.}
Martin Gardner, with a column in \emph{Scientific American}, which was taken up in his book~\cite{Gardner:1991} made the 
surprise test paradox known to a wide audience.
Halpern and Moses~\cite{HalpernMoses:1986} use standard logic, Goedelization and a fixed-point operator
to analyze the situation, and they find that the full, self-referential announcement is contradictory.  
Fitch~\cite{Fitch:1964} makes the Goedelization viewpoint explicit, from a mathematical viewpoint, 
and he comes to the conclusion that the announcement in its full strength must be contradictory.
A further contribution in this direction can be found in~\cite{Kritchman:Raz:2010}, as late as 2010. 
Blinkley~\cite{1968:Binkley:Surprise:Examination} was the first to expressly use modal logic, in the epistemic version, with operators for ``knowing at a certain time'', 
and he reaches the same conclusion as~\cite{Kaplan:Montague:1960}:
Surprise can happen, but not if the students ``know'' beforehand that the announcement is true.
Sorensen~\cite{Sorensen:1982} introduced a version of the paradox that does not mention time,
in terms of modal logic. 
Sorensen devotes a whole chapter in his book ``Blindspots''~\cite{sorensen1988blindspots} to the history of the discussion of the 
story and another chapter to his own analysis in terms of epistemic logic.

Too much has been written and published about this story, so it is not possible to 
give a complete picture of the state of affairs at this point. Let it be enough to say that people from 
mathematics, logic, and philosophy have not found a consensus of what to make of it.
In many papers the conviction is stated that there is a full explanation, which shows there is no paradox, 
only a flawed argument on the side of the students.
The most recent study concerning the story from the point of view of modal logic is~\cite{2025:Baltag:topological:epistemic},
introducing ``topological epistemic logic'' to resolve the ``paradox''.
But the opinions diverge about what exactly the flaw(s) is (are). 
Other works disagree completely and see an actual paradox.

%%%%%%%%%%%%%%%%%%%%%%%%%%%%%%%%%%%%%%%%%%%%%%%%%%%%%%%%%%%%%%%%%%%
\section{Knowledge and surprise}\label{sec:framework}
%%%%%%%%%%%%%%%%%%%%%%%%%%%%%%%%%%%%%%%%%%%%%%%%%%%%%%%%%%%%%%%%%%%

In this section, we set up a framework for describing the events of the week to
provide exact mathematical counterparts for the students ``knowing'' something and for being ``surprised''
of a test on a certain day before the background of this ``knowledge''.

\subsection{Options and reformulation}
%%%%%%%%%%%%%%%%%%%%%%%%%%%%%%%%%%%%%%%%%%%%%%%%%%%%%%%%%%%%%%%%%%%%
The teacher has six options: Give a test on one of the
five days of the week, or give no test at all. 
These options are represented as elements of $R=\{\Mo,\Tu,\We,\Th,\Fr,\none\}$,
called \emph{runs}, and linearly ordered as $\Mo<\dots<\Fr<\none$.
The elements of $D=\{\Mo,\dots,\Fr\}$ are the \emph{days}.
If $d$ is a day, then $d+1$ is the element of $R$ following $d$, 
and $d-1$ is the day before $d$, if it exists.
We identify day $d$ with the run $d$ in which the test is on day $d$. 

The teacher chooses an element $\rT$ of $R$.
It does not matter when this choice is made. 
(Of course the students do not have access to the information what $\rT$ is.)
Some runs may be excluded by ``law'', given as a set $L\subseteq R$,
with the meaning that the teacher can never choose a run from $R-L$.
Set $L$ is always fixed and known to all parties.
We note the following cases.  
\begin{description}
\item{$L=R$: \ }   The standard STP story. 
\item{$L=D$: \ }   The ``\emph{surprising egg paradox}''~\cite{Scriven:1951}.
\item{$L=\{\Mo,\none\}$: \ }   The ``\emph{birthday present paradox}''~\cite{Gardner:1991,Scriven:1951}.
\end{description}
As precise equivalent to ``knowing'' we choose ``being able to prove''. 
This is still informal, since the details of what kinds of proofs are allowed have yet to be spelled out.
The announcement then reads as follows, under the assumption that $L\subseteq R$ has been fixed and is known (to the students, to everybody):
\begin{description}
\item{{\SpL} \ } \begin{minipage}[t]{0.83\textwidth}
My choice $\rT$ is some $d\in D\cap L$. On day $d$ we will have a test, which
will be `surprising' in the sense that on the morning of day $d$
you will not be able to prove that the test is on this day.''
\end{minipage}
\end{description}
Once we have clarified what it means to ``prove something on the morning of day $d$\kern1pt'', 
and with $\rT$ being some fixed value, announcement {\SpL} is either true or false.

\subsection{Formulas, models, p-provability}\label{sec:propositional:formulas}
%%%%%%%%%%%%%%%%%%%%%%%%%%%%%%%%%%%%%%%%%%%%%%%%%%%%
We employ five propositional variables $X_1,\dots,X_5$, where $X_i$ stands for ``there is a test on the $i$-th day of the week''.
Propositional formulas $\varphi$ over $X_1,\dots,X_5$ are built in any standard way. For example, we could define inductively that (i) $\bot$ is a formula,
(ii) $X_i$ is a formula for $1\le i \le 5$, and (iii) for formulas $\varphi$ and $\psi$ also $(\varphi\to\psi)$ is a formula.
With $\neg\varphi:=(\varphi\to\bot)$ and $(\varphi\vee\psi):=(\neg\varphi\to\psi)$
we obtain a standard basis of boolean operators, and $\wedge$ and $\leftrightarrow$ can be defined as usual. Finally we let $\top:=\neg\bot$.
In propositional logic, the models are the assignments $a=(a_1,\dots,a_5)\in\{0,1\}^5$ to the variables. (So in our case there are 32 models.)
We define the truth value $\llbracket\varphi\rrbracket_a$ (in $\{0,1\}$) of a formula $\varphi$ under assignment $a$ in the usual way,%
\footnote{(i) $\llbracket\bot\rrbracket_a=0$, (ii) $\llbracket X_i\rrbracket_{(a_1,\dots,a_5)}=a_i$, $1\le i\le5$, 
(iii) $\llbracket(\varphi\to\psi)\rrbracket_a=\max\{1-\llbracket\varphi\rrbracket_a,\llbracket\psi\rrbracket_a\}$.}
and we write $a\models\varphi$ for $\llbracket\varphi\rrbracket_a=1$.

To express that there cannot be more than one test in the week, we adopt
\[
\bigwedge_{\substack{1\le i<j\le 5}}\!(\neg X_i \vee \neg X_j)\tag*{(Ax$_{\le1}$)}
\]
as an axiom. 
Under (Ax$_{\le 1}$), only six relevant assignments remain, namely those
that have at most one component that is 1,
and these correspond to the elements of $R$.
This correspondence is so natural and obvious that we identify $R$
with the set of these six assignments (i.e., models), which are also called $\Mo,\dots,\Fr,\none$.

We define six basic formulas, which will then be used as auxiliary variables:
\begin{align*}Y_\Mo&:= X_1, \;\; Y_\Tu:= X_2, \;\; Y_\We:= X_3, \;\; Y_\Th:= X_4, \;\; Y_\Fr:= X_5,\\
Y_\none&:= \neg(X_1\vee X_2\vee X_3\vee X_4\vee X_5).
\end{align*}
Clearly, axiom (Ax$_{\le1}$) is equivalent to the following formula, which expresses that exactly one
of the auxiliary variables has value 1:
\[
\bigvee_{r\in R}Y_r\;\; \wedge\; \bigwedge_{\substack{r,s\in R\\r < s}}\!(\neg Y_r \vee \neg Y_s).
\tag*{(Ax$_{=1}$)}
\]
We can now forget about $X_1,\dots,X_5$ and (Ax$_{\le1}$), 
write formulas as boolean combinations of $Y_\Mo,\dots,Y_\Fr,Y_\none$, and use 
(Ax$_{=1}$) as our basic axiom, mostly without mentioning it.
For formulas $\varphi$ and $\psi$ we write
$\varphi\equiv\psi$ for semantic equivalence within the models in $R$, i.e., 
for $\llbracket\varphi\rrbracket_r=\llbracket\psi\rrbracket_r$ for all $r\in R$. 

In order to be able to talk about \emph{provability} in propositional logic (``\emph{p-provability}'') we adopt an arbitrary one of the standard deduction calculi.
If a formula $\varphi$ can be proved from (Ax$_{=1}$), we write $\vdash\!\varphi$.
We will consider arbitrary \emph{axiom systems} $\Ac$, i.e., sets of formulas, always silently assuming that (Ax$_{=1}$) is added to $\Ac$ if necessary, 
and we write $\Ac\vdash\!\varphi$ if $\varphi$ is provable from $\Ac$, for $\varphi$ a formula.
Of course, we will never prove anything in propositional logic. 
Since the deduction calculus is sound and complete, we can argue ``semantically'' on the basis of the six models.
Whatever is the case there [in all models that obey $\Ac$] is also provable [from $\Ac$], and vice versa.
In particular we have $\varphi\equiv\psi$ if and only if $\vdash\,\varphi\leftrightarrow\psi$, for formulas $\varphi$ and $\psi$.
An axiom system $\Ac$ is \emph{inconsistent} if $\Ac\vdash\varphi$ for all formulas $\varphi$,
otherwise $\Ac$ is \emph{consistent}. Clearly $\Ac$ is consistent if and only if there is a model that satisfies all axioms in $\Ac$. 

For each subset $B$ of $R$ we define
\[\chi_B:=\bigvee_{r\in B}Y_r.\] 
With~(Ax$_{=1}$) one can easily show that each formula $\varphi$ is equivalent to some formula $\chi_B$.
Namely, with $R_\varphi:=\{r\in R \mid r\models \varphi\}$ we have the following.
\begin{observation}\label{obs:knowledge:sets:basis}
$\varphi\,\equiv\,\chi_{R_\varphi}$, \ for all formulas $\varphi$.\qed 
\label{obs:standard:form:for:formulas}
\end{observation} 
This observation provides us with standard versions for all formulas.
Actually, we can imagine that the 64 formulas $\chi_B$ are all there are, 
since $\bot\equiv\chi_\emptyset$, $Y_r\equiv\chi_{\{r\}}$, for $r\in R$, and $(\chi_B\to\chi_C)\equiv\chi_{(R-B)\cup C}$,
for arbitrary $B,C\subseteq R$.

In order to translate an assertion S (i.e., a statement with a truth value) into a formula we use \emph{Iverson notation}~\cite{GrahamKnuthPatashnik:1994}:
\[[\,\text{S}\,]=\left\{
\begin{array}{ll}
\top  &  \text{if S is true}\\
\bot  &  \text{if S is false.}
\end{array}
\right.
\]
We introduce suggestive abbreviations for certain formulas. 
For this, we fix a symbol $\cT$ (``choice of the teacher'') and write
\begin{align*}
	\langle \cT\in B\rangle  &\text{\quad for \ \ $\chi_B$}.\\ 
	\intertext{\mbox{We note that $r\models \langle\cT\in B\rangle$ \ $\Leftrightarrow$ \ $r\models \chi_B$ \ $\Leftrightarrow$ \ $r\in B$.
	Similarly, we write}}
	\langle \cT\ge d\rangle  &\text{ \ for \ ${\textstyle\bigvee_{r\ge d}Y_r}$,}
	&
	\langle \cT=r\rangle  &\text{ \ for \ $Y_r$,}
	\\
	\langle \cT\le d\rangle  &\text{ \ for \ ${\textstyle\bigvee_{r\le d}Y_r}$,}
	&
	\langle \cT\neq r\rangle  &\text{ \ for \ $\neg Y_r$.}
\end{align*}

By Observation~\ref{obs:knowledge:sets:basis} we have $\varphi\equiv\langle\cT\in R_\varphi\rangle$, so that
formulas $\langle \cT\in B\rangle$ are sufficient to express everything that can be expressed by formulas in general. 

\subsection{``Knowledge''}\label{sec:knowledge}
%%%%%%%%%%%%%%%%%%%%%%%%%%%%%%%%%%%%%%%%%%%%%%%%%%%%%%%%%%%%%
An axiom system $\Ac$ is a set of formulas containing (Ax$_{=1}$).
What the students ``know'' when they have $\Ac$ is what can be p-proved from $\Ac$. We define:
\begin{definition}\label{def:knowledge:set}
$K\subseteq R$ is an \emph{$\Ac$-p-knowledge set} \  if \ $\Ac\;\vdash\langle\cT\in K\rangle$.
\end{definition}
By what we said before, this definition captures all formulas that are p-provable from $\Ac$, up to equivalence. 
That $K$ is an $\Ac$-p-knowledge set means that if $\Ac$ is satisfied in a model $r$
then $r$ must be in $K$. The structure of the collection of $\Ac$-p-knowledge sets for an axiom system $\Ac$ is very simple.
Because $R$ is always an $\Ac$-p-knowledge set and the collection of $\Ac$-p-knowledge sets is closed under intersection
and taking supersets, 
there is a smallest $\Ac$-p-knowledge set, namely 
\[K_{\Ac}:={\textstyle\bigcap\{K\subseteq R\mid K\text{ is an $\Ac$-p-knowledge set}\}},\]
and $K$ is an $\Ac$-p-knowledge set if and only if $K_\Ac\subseteq K$.
We call axiom systems $\Ac_1$ and $\Ac_2$ \emph{equivalent}, in symbols $\Ac_1\sim\Ac_2$, 
if all formulas in $\Ac_2$ can be p-proved from $\Ac_1$ and the other way round. The following is easily seen.
\begin{observation}\label{obs:axioms:knowledge:sets}For axiom systems $\Ac_1$ and $\Ac_2$ the following are equivalent:
\begin{itemize}
	\item[(i)] $\Ac_1\sim\Ac_2$.
  \item[(ii)] $\Ac_1$ and $\Ac_2$ have the same p-knowledge sets.
  \item[(iii)] $K_{\Ac_1}=K_{\Ac_2}$.\qed
\end{itemize}
\end{observation}
Disregarding equivalences, an axiom system $\Ac$ with $L=K_\Ac$ is called $\Ac_L$. 
One such axiom system is $\{\text{(Ax$_{=1}$)},\langle\cT\in L\rangle\}$.
There are exactly $64$ nonequivalent axiom systems, one for each subset of $R$. 
One is $\Ac_\emptyset\sim\{\text{(Ax$_{=1}$)},\langle\cT\in\emptyset\rangle\}$, the contradictory axiom system;
at the other end there is the weakest axiom system $\Ac_R\sim\{\text{(Ax$_{=1}$)}\}$.

\subsection{Surprise, given an axiom system}\label{sec:def:surprise:L}
%%%%%%%%%%%%%%%%%%%%%%%%%%%%%%%%%%%%%%%%%%%%%%%%%%%%%%%%%%%%%%%%%%%%%%%%%%%%%%%%%%%%%%%
We wish to capture within our framework the situation that a test on day $d$ is ``surprising'', given some axiom system $\Ac$.
The choice $\rT$ of the teacher must satisfy $\Ac$, in other words, we must have $\rT\in K_\Ac$.
If $d\notin K_{\Ac}$, then $d$ cannot be chosen as $\rT$, so it is uninteresting.
So assume $d\in K_{\Ac}$. We ask when a test on day $d$ is ``expected'' relative to $\Ac$, or \emph{$\Ac$-p-expected}.
This is the case if it can be p-proved that the test is on day $d$, 
on the basis of $\Ac$ and the additional information available to the students on the morning of day $d$. 
This information is that it is absolutely impossible that the test is on days $\Mo,\dots,d-1$
(because the students have observed that these days have passed and no test has materialized). 
So the criterion for an $\Ac$-p-expected test is that the following statement can be p-proved from $\Ac$:  
``If no other choices than $d,d+1,\dots,\Fr,\none$ are possible, then there must be a test on day $d$\kern1pt''.
This clearly means the same as that ``no option in $\{d+1,\dots,\Fr,\none\}$ can be chosen by the teacher'' can be p-proved from $\Ac$, which 
is to say that $\{\Mo,\dots,d\}$ is an $\Ac$-p-knowledge set. 

This discussion suggest the following definition.
\begin{definition}\label{def:surprising:run}
$d\in D$ is called \emph{$\Ac$-p-surprising} if $d\in K_\Ac$ and $\Ac\nvdash\langle\cT\le d\rangle$. We let
\[D_{\text{surpr}}^\Ac:=\{d\in D \mid \text{$d$ is \emph{$\Ac$-p-surprising}} \}.\]
\end{definition}
\begin{observation}\label{obs:surprise}
For all $d\in D\cap K_\Ac$ we have: 
\begin{center}
$d$ is $\Ac$-p-surprising \ \ if and only if \ \ $\Ac\cup\{\langle\cT\ge d\rangle\}\;\nvdash\;\langle \cT=d\rangle$.
\end{center}
\end{observation}

\emph{Proof}: Apply the deduction theorem for propositional logic.
\qed 

We note in passing that in Definition \ref{def:surprising:run} any concept of ``time'' has disappeared.
This corresponds to the old observation (see in particular the ``designated student paradox'' from~\cite{Sorensen:1982}) 
that in contrast to first appearances time and the development of ``knowledge'' over time does not play a role in STP.   

Observation~\ref{obs:surprise} expresses nicely that our definition of an $\Ac$-p-surpris\-ing day makes intuitive sense. 
Assume the teacher choses $\rT=d$: the test is on day $d$.
On the morning of day $d$ the students have $\Ac$ and, by observation, that $\rT\ge d$, i.e., $\rT\models \langle\cT\ge d\rangle$.
They can add this information to their ``knowledge'' by adding $\langle\cT\ge d\rangle$ to their axioms,
without overlooking any run, i.e., model, that can possibly be relevant. 
Observation~\ref{obs:surprise} says that choice $d$ is \mbox{$\Ac$-p}-surprising if and only if
one cannot p-prove from this extended axiom system 
that choices $r>d$ are impossible. This corresponds to the intuitive idea of a test on day $d$ being ``surprising'', 
if the information that can be used is given by~$\Ac$.

It will be useful to have a formula that characterizes $\Ac$-p-surprising runs.
At the first glance this seems hard to achieve, since such a formula
must somehow refer to $\Ac$-p-knowledge sets, which involve p-proofs, 
something that looks to be beyond what can be expressed with propositional formulas. 
We employ Iverson notation to force a place for provability statements in propositional formulas and
consider \[[K\text{ is an $\Ac$-p-knowledge set}]=[\Ac\vdash\langle \cT\in K\rangle]=[K_\Ac\subseteq K],\] 
which equals $\top$ if $K$ is an $\Ac$-p-knowledge set and equals $\bot$ if it is not.  
With this, we define:%
\footnote{Stretching notation in a harmless way, we take $\emptyset-\{\max \emptyset\}$ to mean $\emptyset$.}
\begin{definition}\label{def:surprising:run:formula}\quad\\[-12pt]
\[%\label{eq:def:D:surpr}
\sigma^\Ac:=\;\; \bigwedge_{K\subseteq R}\Bigl([K_\Ac\subseteq K]\to\langle\cT\in K{-}\{\max K\}\rangle\Bigr).
\]
\end{definition}
\begin{lemma}\label{lem:formula:surprise}
For $d\in D$ the following are equivalent: 
\begin{itemize}
	\item[(i)]  $d\models \sigma^\Ac$.
	\item[(ii)]	 For all $\Ac$-p-knowledge sets $K$: $d\in K$ and $d\neq \max K$.
	\item[(iii)] $d\in K_\Ac$ and $d\neq \max K_\Ac$.
	\item[(iv)] $d\in D_{\text{surpr}}^\Ac$. 
\end{itemize}
\end{lemma}
\emph{Proof}: That (i) and (ii) are equivalent follows from the definition of $\sigma^\Ac$:
Factor $\pi^\Ac_K:=([K_\Ac\subseteq K]\to\langle\cT\in K{-}\{\max K\}\rangle)$ in the conjunction in $\sigma^\Ac$
acts like a ``conditional formula'': If $K$ is not an $\Ac$-p-knowledge set, $\pi^\Ac_K\equiv\top$,
thus it is nonexistent in the context of the conjunction, otherwise $\pi^\Ac_K\equiv\langle\cT\in K{-}\{\max K\}\rangle$.
\\[2pt]
``(ii) $\Rightarrow$ (iii)'': Assume $d\in K-\{ \max K\}$ for all $\Ac$-p-knowledge sets $K$. Then $d\in K_\Ac$,
and $d\neq \max K_\Ac$, because $K_\Ac$ is itself an $\Ac$-p-knowledge set.
\\[2pt]
``(iii) $\Rightarrow$ (ii)'': Assume $d\in K_\Ac-\{ \max K_\Ac\}$.
Then $d\in K$ for all $\Ac$-p-knowledge sets $K$. And it is impossible that $d=\max K$ for some $\Ac$-p-knowledge set $K$, 
since by $K\supseteq K_\Ac$ this would imply $d\ge \max K_\Ac$, hence $d=\max K_\Ac$.
\\[2pt]
``(iii) $\Rightarrow$ (iv)'': Assume $d\in K_\Ac$ and $d\neq \max K_\Ac$. 
Then $\{\Mo,\dots,d\}$ cannot be an $\Ac$-p-knowledge set, since otherwise we would get $d=\max K_\Ac$. 
This means that $\Ac\nvdash\langle\cT\le d\rangle$, i.e., $d\in D_{\text{surpr}}^\Ac$.
\\[2pt]
``(iv) $\Rightarrow$ (iii)'': Assume $d\in D_{\text{surpr}}^\Ac$. Then $d\in K_\Ac$ and $\Ac\nvdash\langle\cT\le d\rangle$, 
which means that $K_\Ac\not\subseteq\{\Mo,\dots,d\}$, i.e., $\max K_\Ac > d$.    
\qed

``Surprise'' becomes impossible once there is an $\Ac$-p-knowledge set of size 1. 

\begin{observation}\label{obs:fringe}
If $\Ac\vdash\langle \cT=r\rangle$ for some $r\in R$ then $D_{\text{surpr}}^\Ac=\emptyset$.
\end{observation}
\emph{Proof}: Assume $\Ac\vdash\langle \cT=r\rangle$. This implies $K_{\Ac}\subseteq\{r\}$.
If $K_{\Ac}=\emptyset$, we are done, and if $K_{\Ac}=\{r\}$, we have
$D_{\text{surpr}}^\Ac=D\cap K_{\Ac}-\{\max K_{\Ac}\}=\emptyset$.
\qed

This finishes the description of STP in terms of propositional logic,
given an axiom system $\Ac$. We still have to find out what $\Ac$ should be in view of {\SpL}.
In Section~\ref{sec:students} we discuss the (incremental) strategy the students use to gather information about $\Ac$, 
and in Section~\ref{sec:general:analysis} we provide a general analysis.

%%%%%%%%%%%%%%%%%%%%%%%%%%%%%%%%%%%%%%%%%%%%%%%%%%%%%%%%%%%%%%%%%%%
\section{The argument of the students}\label{sec:students}
%%%%%%%%%%%%%%%%%%%%%%%%%%%%%%%%%%%%%%%%%%%%%%%%%%%%%%%%%%%%%%%%%%%%

We now look at the argument of the students in detail, show how it 
can be understood as establishing tighter and tighter conditions on the 
axiom system to be used, and expose all its flaws.  

\subsection{Self-referentiality of the announcement}\label{sec:announcement:self}
%%%%%%%%%%%%%%%%%%%%%%%%%%%%%%%%%%%%%%%%%%%%%%%%%%%%%%%%%%%%%%%%%%%%%%%%
We return to the original version of STP, with $L=R$.
Recall announcement {\SpL} (with $L=R$). For ``surprising'' we would like to use ``$\Ac$-p-surprising''
for a suitable axiom system $\Ac$. But what should this axiom system be?
The story seems to tell us that the truth of {\SpL} with $\rT$ gives us 
a basis for finding new p-knowledge sets, i.e., to show that certain formulas must be p-provable from $\Ac$.
Of course, then, more days may turn out to be not $\Ac$-p-surprising,
which means more information about $\Ac$, and so on. The following version of the announcement
expresses that whatever can be proved from the announcement must be reflected directly in the axiom system $\Ac$,
so that it influences what is considered $\Ac$-p-surprising. 

\begin{description}
\item{{\Spiter} \ } \begin{minipage}[t]{0.82\textwidth}
``My choice $\rT$ is a day $d$ in $D_{\text{surpr}}^\Ac$, where $\Ac$ fulfills the following:\\[2pt] 
Whenever for a set $B\subseteq R$ you can prove that 
 $r\models\Ac$ entails $r\in B$, for all $r\in R$,
then $\Ac\vdash\langle \cT\in B\rangle$. 
In your proofs you may use that $\rT\in D_{\text{surpr}}^\Ac$ is true.''
\end{minipage}
\end{description}
 
Here ``prove'' and ``entails'' refers to proofs in the usual mathematical sense.
Note that this does not open as vast a field as it looks at the first glance, since the mathematical
domain that we cover is only axiom systems, proofs, models, and validity
for the tiny part of propositional logic cut out by variables $Y_r$, $r\in R$, and (Ax$_{=1}$).
Using that $r$ does not appear in ``$\rT\in D_{\text{surpr}}^\Ac$'',  the condition on $\Ac$ in {\Spiter} is easily seen to be equivalent to
\begin{equation}\label{eq:condition}
\text{if \ $\rT\in D_{\text{surpr}}^\Ac$ \ $\Rightarrow$ \ $\Ac\vdash\langle \cT\in B\rangle$, \ \ then \ \ $\Ac\vdash\langle \cT\in B\rangle$\kern1pt.}
\end{equation}
Condition~\eqref{eq:condition} looks like a closure condition on $\Ac$.
We shall see in due course what its effects are.

\subsection{Axioms ``established'' by the students}\label{sec:axiom:system:students}
%%%%%%%%%%%%%%%%%%%%%%%%%%%%%%%%%%%%%%%%%%%%%%%%%%%%%%%%%%%%%%%%%%%%%%%
We translate the argument of the students in the original story into our framework.
The students take {\Spiter} as a specification of some (as yet undetermined) axiom system $\Ac$.
They assume the teacher has chosen some $\rT$ (equally unknown)
from $D_{\text{surpr}}^\Ac$.
By Lemma~\ref{lem:formula:surprise} we can write $\rT\in D_{\text{surpr}}^\Ac$ as
\begin{equation}\label{eq:hypothesis}
\rT\in D\text{ \ and \ }\rT\models\sigma^\Ac.
\end{equation}
The students use \eqref{eq:condition} and~\eqref{eq:hypothesis} to iteratively
derive information about $\Ac$. In detail, this proceeds as follows. 
\begin{enumerate}\setcounter{enumi}{-1}
\item {(Assumptions.)} \ Axiom system $\Ac$ fulfills \eqref{eq:condition} and~\eqref{eq:hypothesis}.
\item (The ``initial step''.) \  The students say: ``From~\eqref{eq:hypothesis} we get $\rT\in D$.
\begin{itemize}
\item[(1a)] This means that the teacher cannot choose {\none} as her $\rT$.''
\end{itemize}
From this they go on to
\begin{itemize}
\item[(1b)] ``Under $\Ac$, choice {\none} is impossible.''
\end{itemize}
This they read as $\Ac\vdash\langle \cT\le\Fr\rangle$, 
which by~\eqref{eq:condition} they consider as a property of $\Ac$ newly established.
\item (The ``no-Friday step''.) \ For $d=\Fr,\Th,\We,\Tu$, in this order, the students argue in the manner of a proof by induction.%
\footnote{Readers familiar with the concept should note that this is not backward induction, 
but a plain, direct induction, ``proving'' more and more properties of $\Ac$.}
	From the previous step they (believe that they) have $\Ac\vdash\langle \cT\le d\rangle$, i.e., $\{\Mo,\dots,d\}$ is an $\Ac$-p-knowledge set. 
	They say: ``Together with $\rT\models\sigma^\Ac$, which holds by~\eqref{eq:hypothesis}, and looking at Definition~\ref{def:surprising:run:formula} 
	we get that $\rT\models\langle\cT\neq d\rangle$, i.e., $\rT\neq d$. 
\begin{itemize}
\item[(2a)] This means that the teacher cannot choose $d$ as her $\rT$.''
\end{itemize}
From this they go on to
\begin{itemize}
\item[(2b)] ``Under $\Ac$, the teacher cannot choose $d$ at all.''
\end{itemize}
This they read as $\Ac\vdash\langle \cT\ne d\rangle$,  
which by~\eqref{eq:condition} they consider as a property of $\Ac$ newly established.
Combining this with $\Ac\vdash\langle \cT\le d\rangle$ from the previous round they obtain $\Ac\vdash\langle \cT\le d-1\rangle$.
\item (The ``final step''.) The result of Step 2.{} for $d=\Tu$ is $\Ac\vdash\langle\cT=\Mo\rangle$.
Now Observation~\ref{obs:fringe} gives that this entails $D_{\text{surpr}}^\Ac=\emptyset$,
which contradicts the initial assumption $\rT\in D_{\text{surpr}}^\Ac$.%
\footnote{We avoid the last big logical mistake of the students, which leads them to 
believe they had proved that there cannot be a test.}
\end{enumerate}

Isolating the central steps in the argument of the students,
and transforming them into formulas, we see that 
their argument boils down to the idea that  
$\rT\in D_{\text{surpr}}^\Ac$ entails that the five formulas 
\begin{equation}\label{eq:students:consequences}
\langle \cT\in D\rangle\text{ \ and \ }[\,\Ac\vdash\langle\cT\le d\rangle\,]\to\langle\cT\ne d\rangle\text{, \ for }d\in\{\Tu,\dots,\Fr\}\text{,}
\end{equation}
are all consequences of $\Ac$.
Our discussion above can easily be turned into a proof of the following. 
\begin{proposition}\label{prop:students:overall}
Let $\Ac$ be an axiom system from which the five formulas in~\eqref{eq:students:consequences}
are provable. 
Then $D_{\text{surpr}}^\Ac=\emptyset$, i.e., no $d\in D$ satisfies $d\models \sigma^\Ac$.\qed
\end{proposition}
\begin{corollary}\label{cor:sigma:inconsistent}
If axiom system $\Ac$ implies $\langle \cT\in D\rangle$ and $\sigma^\Ac$, then $\Ac$ is inconsistent.
\end{corollary}
\emph{Proof}: Looking at Definition~\ref{def:surprising:run:formula} one easily sees that
all formulas in~\eqref{eq:students:consequences} are provable from $\langle \cT\in D\rangle\wedge\sigma^\Ac$.\qed 

It is a fact that the approach of demanding that 
$\langle \cT\in D\rangle$ and $\sigma^\Ac$ must be axioms in $\Ac$ 
(or something equivalent to it) is taken in many of the previous works
on STP, 
e.g.~\cite{Fitch:1964,HalpernMoses:1986,Kaplan:Montague:1960,Kritchman:Raz:2010}.
Many authors have tried to escape the resulting inconsistency by taking refuge to weaker logics,%
\footnote{Beware: {\Spiter} does not only use propositional logic, but standard mathematical proofs on top of it.}
like different versions of modal logic, in particular epistemic logic of some sort~\cite{Aldini:et:al:2023,2008:Baltag:epistemic,1968:Binkley:Surprise:Examination,%
Gerbrandy:2007,Hall:1999,2018:Holliday:epistemic,Murzi:2021,2006:VanDithmarsch:success}. 
This escape route is also mentioned in Chow's recipe cited in Section~\ref{sec:story:paradox}.
We choose a different path, in sticking to full propositional logic. 
We just try axioms weaker than those in~\eqref{eq:students:consequences}. 
Of course, it must be carefully argued that this still results in a faithful translation of the story into mathematical terms.
This will be done in Section~\ref{sec:general:analysis}. Here, we proceed in a more traditional way and show that
$\rT\in D_{\text{surpr}}^\Ac$ being true does not imply any of the 
formulas in~\eqref{eq:students:consequences}. 

This means that the students are very much in error, and that they make five unjustifiable steps.

\subsection{Crushing the argument of the students}\label{sec:crush}
%%%%%%%%%%%%%%%%%%%%%%%%%%%%%%%%%%%%%%%%%%%%%%%%%%%%%%%%%%%%%%%
The following theorem summarizes our case against the argumentation of the students,
as elaborated in Section~\ref{sec:axiom:system:students}. (Presumably the readers 
can identify the dubious steps in that argumentation on their own.)   
\begin{theorem}\label{thm:students:are:wrong}
Axiom system $\Ac_R=\{\text{(Ax$_{=1}$)}\}$ has the following properties.
\begin{itemize}
	\item[(a)] $r\models\Ac_R$ for all $r\in R$. Hence $K_{\Ac_R}=R$, and $\Ac_R$ is consistent. 	
	\item[(b)] $D_{\text{surpr}}^{\Ac_R}=D$.
	\item[(c)] $\Ac_R\nvdash\langle\cT\in D\rangle$.
\end{itemize}
Further we have:
\begin{itemize}
	\item[(d)] For $d\in D-\{\Mo\}$ axiom system $\Ac_d:=\Ac_{\{\Mo,\dots,d\}}=\{\text{(Ax$_{=1}$)},\langle \cT\le d\rangle\}$ 
	has the following properties: 
	\begin{itemize}
		\item[(i)] $r\models \Ac_d$ for all $r\le d$.\\ Hence $K_{\Ac_d}=\{\Mo,\dots,d\}$, and $\Ac_d$ is consistent;
	\item[(ii)] $D_{\text{surpr}}^{\Ac_d}=\{\Mo,\dots,d-1\}$;  
	\item[(iii)]$\Ac_d\;\;\nvdash\;\;[\,\Ac_d\vdash\langle\cT\le d\rangle\,]\to\langle\cT\ne d\rangle$.  
	\end{itemize}
\end{itemize}
\end{theorem} 
\emph{Proof}:
(a) That $r\models\Ac_R$ for all $r\in R$ is trivial.
If an axiom system has a model, it is consistent. \\[2pt]
(b) Apply Definition~\ref{def:surprising:run} with (a).\\[2pt]
(c) $\Ac_R\vdash\langle \cT\in D\rangle$ would mean that $d$ is an $\Ac_R$-p-knowledge set, which it is not, by (a).\\[2pt]
(d) Let $d\in\{\Tu,\We,\Th,\Fr\}$. With $\Ac_d=\{\text{(Ax$_{=1}$)},\langle \cT\le d\rangle\}$ we have:\\[2pt] 
(i) This is trivially true.\\[2pt] (ii) Apply
Definition~\ref{def:surprising:run} with (i).\\[2pt] 
(iii) We have $\Ac_d\vdash\langle\cT\le d\rangle$, 
but $d\not\models\langle\cT\ne d\rangle$. Thus $d\not\models[\,\Ac_d\vdash\langle\cT\le d\rangle\,]\to\langle\cT\ne d\rangle$, but $d\models\Ac_d$ by (i). 
Hence $\Ac_d\;\nvdash\;[\,\Ac_d\vdash\langle\cT\le d\rangle\,]\to\langle\cT\ne d\rangle$.
\qed
\begin{remark} 
In part (d), the case $d=\Mo$ is conspicuously missing. 
Although the students in the story say ``Just like for Friday, we can exclude Monday'', this case is different. 
If $\Ac\vdash\langle \cT=\Mo\rangle$ for the axiom system $\Ac$ in question,
Observation~\ref{obs:fringe} says that $\rT\in D^\Ac_{\text{surpr}}$ cannot be true with any $\rT$,
so we have reached a contradiction already here. This means that the whole situation is contradictory, 
i.e., nonexistent, and there everything is provable, including $\Ac\vdash\langle\cT\neq\Mo\rangle$.
So this last conclusion looks strange, but it is not erroneous.
\end{remark}
Theorem~\ref{thm:students:are:wrong} and its proof are written in a very compact way, so
some more explanations are warranted.
We first point out that (a)--(c) show that with $\Ac=\Ac_R$ an arbitrary value $\rT\in D$ will
realize {\Spiter} with $\Ac$, but $\Ac\nvdash\langle\cT\in D\rangle$, so that the initial step of the students is not sound.%
\footnote{For the earliest paper that claims to show this, in a less formalized manner, see~\cite{Quine:1953}.} 
Further, (d)(i)--(iii) show, for arbitrary $d\in\{\Tu,\dots,\Fr\}$,
that with $\Ac=\Ac_d$ an arbitrary value $\rT\in \{\Mo,\dots,d-1\}$ will
realize {\Spiter} with $\Ac$, but $\Ac\nvdash\langle\cT\neq d\rangle$, so that the ``no-Friday'' step of the students for $d$ is not sound either. 

The key to understanding how the argument of the students creates a contradiction where
intuitively there should not be any is to observe carefully the difference between p-provability statements regarding $\Ac$
and provability statements that involve~$\rT$.

\textbf{The initial step:} \ The assumption $\rT\in D_{\text{surpr}}^\Ac$ is correctly weakened to $\rT\in D$,
with correct interpretation as (1a). This refers to \emph{one} model of $\Ac$, namely $\rT$.  
Now, in a wild jump, the students claim (1b), which says 
that $r\in D$ must be true for \emph{all} models $r$ of $\Ac$,
which is then correctly translated back into technical language as $\Ac\vdash\langle \cT\le\Fr\rangle$.
But this property of the (unknown) axiom system $\Ac$ is not a consequence of 
$\rT\in D$, as demonstrated in part (c) of the theorem. --
To put it in everyday words: It does not hurt if model $r=\none$ is also admitted by $\Ac$, 
since the teacher can still choose a ``good'' run $\rT$, i.e., a day $d$.

\textbf{The ``no-Friday step'':} \ Let us assume, with the students, that $\rT\in D_{\text{surpr}}^\Ac$ and that $\Ac\vdash\langle\cT\le d\rangle$, for some day $d>\Mo$.
The students correctly conclude $\rT\neq d$, and in (2a) interpret this correctly 
as saying that the teacher cannot chose $d$ as her $\rT$. 
In a similar jump as in the initial step they then conclude (2b), which says that the teacher cannot choose $d$ at all under axiom system $\Ac$,
and correctly translate this back as $\Ac\vdash\langle\cT\neq d\rangle$.
But part (d) of the theorem shows that this is not a consequence of $\rT\neq d$. 
Similarly as for the initial step we can say that it does not hurt
if model $d$ is around: The teacher can still choose some $\rT<d$ and make $\rT\in D_{\text{surpr}}^\Ac$ true. 
(Note that the last sentence is wrong for $d=\Mo$, which is why 
this case is not included in the ``no-Friday step''.)

In somewhat more leisurely wording, we can summarize the errors of the students
as follows (disregarding the obvious logical mistake at the end): 
The sentence ``Run $\bar r$ \emph{cannot be chosen} (in the given circumstances)'' 
can have two different interpretations. One is that $\rT\neq \bar r$ must hold,
and the other is that $\Ac\vdash\langle\cT\neq \bar r\rangle$ must hold. 
The students of the story use everyday language to discuss the situation, 
which makes it possible to blur the difference between the two and thus 
create the illusion one could draw certain conclusions about $\Ac$ that in truth are mathematically unsound. 

This finishes the analysis of the argument of the students.
They are wrong in their overall conclusion, and they make five well-camouflaged unjustified steps underway and 
a glaring logical blunder at the end. (This multitude of errors on the side of the students may be one reason for
so many diverging ``solutions'' of the ``paradox'' having been proposed.)
The crux of the matter is that the impreciseness of everyday language
causes the students (and us, the listeners) to see a chain of inferences where actually every single link is broken. 

\begin{remark}\label{rem:arbitrary:Ac} 
Finally, we comment on which axiom systems $\Ac$ satisfy {\Spiter} with some $\rT$. 
The 64 non-equivalent axiom systems split into two classes.
\begin{description}
\item{Case 1:} $|K_\Ac|\ge2$. -- In this case let us assume that the teacher chooses some $\rT\in K_\Ac-\{\max K_\Ac\}$. 
Then $\rT\in D_{\text{surpr}}^\Ac$ is true. This means that the hypothesis ``$\rT\in D_{\text{surpr}}^\Ac$ is true''
in {\Spiter} is worthless, and~\eqref{eq:condition} reduces to the following: 
``Whenever for a set $B\subseteq R$ you can prove that $r\models\Ac$ entails $r\in B$, for all $r\in R$,
then $\Ac\vdash\langle \cT\in B\rangle$.'' But this is a triviality, and it does not 
imply any further restrictions on $\Ac$. 
Since it is conceivable that the teacher acts in this way, the students cannot exclude this behavior, 
and they cannot infer anything further about $\Ac$ from {\Spiter}.
\item{Case 2:} $|K_\Ac|\le1$. -- By Observation~\ref{obs:fringe} this entails $D_{\text{surpr}}^\Ac=\emptyset$,
and the announcement ``$\rT\in D_{\text{surpr}}^\Ac$'' in {\Spiter} is wrong.
No $\rT$ can realize {\Spiter}, for trivial reasons: 
Either the teacher cannot choose anything at all, 
or she has only one choice, and this leads to ``no test'' or to an expected test. 
The only situations in which {\Spiter} makes it possible to prove something 
are these uninteresting cases! 
\end{description} 
\end{remark}

%%%%%%%%%%%%%%%%%%%%%%%%%%%%%%%%%%%%%%%%%%%%%%%%%%%%%%%%%%%%%%%%%%%
\section{A general analysis}\label{sec:general:analysis}
%%%%%%%%%%%%%%%%%%%%%%%%%%%%%%%%%%%%%%%%%%%%%%%%%%%%%%%%%%%%%%%%%%%
With this section we return to the general situation set up in Section~\ref{sec:def:surprise:L},
with an arbitrary law set $L\subseteq R$,
and we aim at characterizing the axiom systems described by {\SpL}.
We have seen so far that the argument of the students in the story is flawed.
Could there be other arguments that still lead to a contradiction?
To answer this, we could expand the ideas of Remark~\ref{rem:arbitrary:Ac}
into a proof that the students can never establish more formulas $\langle\cT\in B\rangle$
than those that follow from $\Ac_L$ because $L\subseteq B$.
Here we choose a different route and undertake a general analysis, which shows that
no proof whatsoever can give the students a tighter axiom system than what they start from,
excepting in the special and quite uninteresting case where the announcement {\Sp} is false for all $\rT$.  
This will provide an alternative and instructive view of the situation in STP.

\subsection{The self-referring announcement}\label{sec:self:referring}
%%%%%%%%%%%%%%%%%%%%%%%%%%%%%%%%%%%%%%%%%%%%%%%%%%%%%%%%%%%%%%%%%%%%%%%
We rephrase the announcement once again, this time emphasizing that the students must not ``invent'' axioms as they like, 
but must prove that they follow from {\SpL}.
\begin{description}
\item{{\SpsrL} \ } \begin{minipage}[t]{0.86\textwidth}
``My choice $\rT$ is a day $d$ in $D_{\text{surpr}}^\Ac$, where
$\Ac$ satisfies the following:
$\Ac\vdash\langle\cT\in B\rangle$ if and only if from $\Ac_L\subseteq\Ac$ and  
{\SpsrL} being true it can be proved that only runs in $B$ can be chosen.''
\end{minipage}
\end{description} 
As before, ``it can be proved'' refers to arbitrary mathematical proofs. This formulation clearly exposes the 
mutual influence, or ``self-reference'', in the announcement of the teacher: 
It uses $\Ac$, and at the same time it is an implicit specification of $\Ac$.%
\footnote{Kaplan and Montague~\cite{Kaplan:Montague:1960} 
use a fixed-point operator in modal logic for defining a corresponding formula.}
By Observation~\ref{obs:axioms:knowledge:sets} 
there are 64 nonequivalent axiom systems $\Ac$.
We will have to find out which of them satisfy {\SpsrL}.

\subsection{One-move games and rational teachers}\label{sec:one:move}
%%%%%%%%%%%%%%%%%%%%%%%%%%%%%%%%%%%%%%%%%%%%%%%%%%%%%%%%%%%%%%%%%%%%%%%%
Like all formulations discussed before announcement {\SpsrL} 
states that the (unknown) run $\rT$ chosen by the teacher has some property. 
For gathering intuition about situations involving such announcements we look at simple games.  

\begin{example}[(Checkmate in one move)] Imagine a chess position where White moves next and White can checkmate Black in one move. 
The White player could tell the Black player: ``I will checkmate you in one move.''
Normally, the White player will have a choice between different moves,
and the standard interpretation of such an announcement is that 
she wants to say that at least one of the choices results in a checkmate. 
In this case one will regard her announcement as ``true''. It does not cause problems
if there are other, non-winning moves. 
An interesting aspect is also that if White can checkmate Black in one move 
no announcement is needed at all to inform the Black player of this, 
since he can find it out for himself from the position on the board. 
Finally, we note that \emph{having a winning move} will be the same as \emph{making a winning move}, 
if the White player behaves ``rationally'' and will choose a winning move if there is one. 
\end{example}
 
We return to the teacher who chooses some run $\rT$ from $L$.  
Without asking where it comes from we assume there is a subset $K_0$ of $L$ 
such that she is allowed to choose an arbitrary run in $K_0$, but the runs in $R-K_0$ are forbidden.
We imagine that the teacher has some objective $\Oc$, which is so that 
some $r$ from $R$ realize $\Oc$ and the others do not. (``Some'' includes ``possibly none, possibly all''.)
Let $S_\Oc=\{r\in R\mid r\in K_0\text{ and $r$ realizes }\Oc\}$ be the set of ``successful'' choices.%
\footnote{In our story the goal of the teacher is to make {\SpsrL} true,
hence $S_\Oc$ is $D_{\text{surpr}}^\Ac$ for an axiom system $\Ac\supseteq\Ac_L$ that satisfies the specification in {\SpsrL}.}
Then a \emph{rational teacher} is one who always chooses a run in $S_\Oc$ if this set is nonempty.%
\footnote{Readers should be aware that being a \emph{rational teacher} is not about a frame of mind, 
but that this label only describes a rule that governs the teacher's behavior for certain situations.
Obviously, the term is chosen so as to evoke the concept of \emph{rational agents} in economics or game theory,
who will always make choices that are most advantageous for them.
In the book~\cite{Gintis:2014}, H.~Gintis explores in great depth (knowledge about) rational behavior in games.}
Whether or not there are runs in $K_0-S_\Oc$ is irrelevant. 
The following is immediate from this definition.
\begin{observation}[(Rationality principle)]
For a rational teacher we have in the situation just described:
\begin{center}
\colorbox{shadethmcolor}{The teacher chooses a run from $S_\Oc$ \ \ if and only if \ \ $S_\Oc\neq\emptyset$.}
\end{center}
(The direction ``$\Rightarrow$'' is true for arbitrary teachers.) 
\end{observation}
It is said nowhere in our story whether the teacher behaves rationally or not.
We (and the students) are well advised to assume a rational teacher. 
We also assume that the teacher is telling the truth, i.e., that {\SpL} with $\rT$ is true. 
Now Corollary~\ref{cor:sigma:inconsistent} tells us that the idea to regard the validity of $\sigma^\Ac$, 
which expresses rational behavior, as a mathematical postulate that can be used in proofs
(also see~\cite{Sober:1998}) leads straightaway into the well-known contradiction.
What should we do in this situation?

Let us brief\nolig ly return to the chess example from above. 
It is reasonable to assume that the White player is rational, and that she makes good on her announcement, 
but intuitively there is no way to prove it as soon as alternative non-winning moves are around. 
Claiming that it can be considered a given fact that she makes a winning move is unrealistic and contradicts common sense. 
In the case of the STP the same is true, but it is even worse: Assuming that a surprising test is a fact to be used in proofs leads into a contradiction.
This amounts to a proof that {\SpsrL} cannot be true if rational behavior is such a given fact.
Of course, it is not astonishing at all that forcing the teacher to behave in a certain way makes her actions predictable
and so surprise becomes impossible. We have the choice: Insist that rational behavior is a fact to be used in proofs
or insist that with {\SpsrL} the teacher speaks the truth. It seems reasonable to opt for the latter, 
since it is intuitive and since the story explicitly mentions that the teacher does not lie,
while it does not say anything about rational behavior.

\subsection{Characterizing surprising days}\label{sec:successful:choices}
%%%%%%%%%%%%%%%%%%%%%%%%%%%%%%%%%%%%%%%%%%%%%%%%%%%%%%%%%%%%%%%%%%%%%%%%%%%%%
Assume an arbitrary law set $L\subseteq R$ is given. 
A day $d$ is called \emph{$L$-surprising} if $d\in D_{\text{surpr}}^\Ac$, 
for an axiom system $\Ac\supseteq \Ac_L$ as specified in {\SpsrL}.
In this section we will see that this $\Ac$ exists and is uniquely determined, 
that we can find what it is and then
read off the set of $L$-surprising days. 

If we assume that the teacher acts rationally, we can replace 
``My choice $\rT$ is a day $d$ in $D_{\text{surpr}}^\Ac$, \dots{}'' in {\SpsrL} by ``There is some day $d$ in $D_{\text{surpr}}^\Ac$, \dots{}'',
or ``$D_{\text{surpr}}^\Ac\neq\emptyset$''.
This leads to the following (new) ``existential'' version of the announcement: 
 
\begin{description}
\item{{\SpsreL} \ } \begin{minipage}[t]{0.82\textwidth}
``$D_{\text{surpr}}^\Ac\neq\emptyset$, where axiom system $\Ac$ is determined by the following condition:
$\Ac\vdash\langle\cT\in B\rangle$ if and only if from $\Ac_L\subseteq\Ac$ and  
{\SpsreL} being true it can be proved that only runs in $B$ can be chosen.''
\end{minipage}
\end{description}

Note that this formulation is at most as strong as {\SpsrL}. 
Since we (with the students) must take into account that the teacher might behave rationally 
we cannot read off more than {\SpsreL} from announcement {\SpsrL}.  

This reformulation changes everything:
The announcement turns into a statement that is either true or false, 
and all the difficulties with the unknown $\rT$ disappear.
It only remains to give the phrase ``prove something from $\Ac_L\subseteq\Ac$ and  
{\SpsreL} being true'' a precise meaning. This we will do next. 

Recall Definition~\ref{def:surprising:run:formula}. The formula $\langle\cT\in D\rangle\wedge\sigma^\Ac$, i.e., 
\[\langle\cT\in D\rangle\;\;\wedge\bigwedge_{K\subseteq R}\!\Bigl([K\text{ is an $\Ac$-p-knowledge set}]\to\langle\cT\in K{-}\{\max K\}\rangle\Bigr)\]
characterizes the elements of $D_{\text{surpr}}^\Ac$. We replace the symbol $\cT$ by $d$, to obtain%
\footnote{The $\cT$ in ``$[\Ac\vdash\langle\cT\in K\rangle]$'', for ``$[K\text{ is  $\Ac$-p-knowledge set}]$'' 
is ``protected'' by the Iverson brackets and the ``$\Ac\vdash\dots$'', and is not replaced, of course.}
\[[ d\in D]\;\;\wedge\bigwedge_{K\subseteq R}\!\Bigl([\Ac\vdash\langle\cT\in K\rangle]\to[ d\in K{-}\{\max K\}]\Bigr),\]
and then define, emulating the existential quantifier in ``$\exists d\in D_{\text{surpr}}^\Ac\dots$'' by a disjunction:
\[\tau^\Ac\;\; := \;\;\bigvee_{d\in D}\;\bigwedge_{K \subseteq R}\!\Bigl([\Ac\vdash\langle\cT\in K\rangle]\to[ d\in K{-}\{\max K\}]\Bigr).\]
(The part $[d\in D]$ equals $\top$ for $d\in D$ and can be dropped.) Since Iverson brackets are $\top$ or $\bot$,
formula $\tau^\Ac$ evaluates to $\top$ or $\bot$, depending on the axiom system $\Ac$.

\begin{proposition}\label{prop:tau:basic}
For an arbitrary axiom system $\Ac\supseteq\Ac_L$ we have: 
\begin{center}
$D_{\text{surpr}}^\Ac\neq\emptyset$ \ \ \ $\Leftrightarrow$ \ \ \ $\tau^\Ac\equiv\top$.
\end{center}
\end{proposition}
\emph{Proof}: Straightforward, with Lemma~\ref{lem:formula:surprise}.\qed

At last, we can give a fully precise version of announcement {\SpsreL}. 
It says the same as 
\begin{description}
\item{{\SpastL} \ \ \ } \begin{minipage}[t]{0.9\textwidth}
``$D_{\text{surpr}}^\Ac\neq\emptyset$, where the
axiom system $\Ac$ satisfies $\Ac\sim\Ac_L\cup\{\tau^\Ac\}$.''
\end{minipage}
\end{description}
If this announcement is true, a rational teacher will choose an element of  $D_{\text{surpr}}^\Ac$. 
In this compact form, the self-referential structure of the announcement becomes crystal clear.
But because $\tau^\Ac$ can only be $\top$ or $\bot$, it is equally clear that 
the self-reference has a very simple structure
and that it is quite easy to give an explicit description of all the possibilities. 

\begin{theorem}\label{thm:rational} Let $L\subseteq R$ be given. Then we have:  
\begin{itemize}
	\item[(a)] If $D\cap L-\{\max L\}\neq\emptyset$
then all consistent axiom systems $\Ac$ as described in {\SpastL} are equivalent to $\Ac_L$,
and $D_{\text{surpr}}^\Ac = D\cap L-\{\max L\}$.
	\item[(b)] If $D\cap L\subseteq\{\max L\}$ then all axiom systems $\Ac$ that satisfy the specification in {\SpastL} are inconsistent, 
	and $D_{\text{surpr}}^\Ac =\emptyset$.
\end{itemize}
\end{theorem}
\emph{Interpretation}: (a) If $D\cap L-\{\max L\}\neq\emptyset$, 
and if the teacher, behaving rationally, chooses a day in $D\cap L-\{\max L\}$ as her $\rT$,
then with {\SpsrL} she has spoken the truth. The students have a test and they are surprised on the day of the test.
This version of surprise has an immediate intuitive meaning:
Part (a) says that the students cannot prove that there are any restrictions (axiom systems, ``knowledge'') narrower than the law set $L$.
Since for $d\in D\cap L$ with $d<\max L$ one cannot prove 
$\Ac_L\cup\{\langle\cT\ge d\rangle\}\vdash\langle\cT\neq\max L\rangle$,
on the morning of day $d$ the students cannot exclude that the teacher chooses $\max L$ (even though this choice does not satisfy {\SpsrL}).
\\[2pt]
(b) If $D\cap L\subseteq\{\max L\}$ then no $r$ can make {\SpsrL} true with an axiom system $\Ac$ that fits the specification.
(If $L=\{\none\}$, the only possible choice is $\rT=\none$. If $L=\{d\}$ for a day, then the only possible choice is $\rT=d$, 
and the students can prove this, so they are not surprised when the test comes on day $d$.)

\emph{Proof}: Let $L$ be given, and assume that $\Ac\sim\Ac_L\cup\{\tau^\Ac\}$.
By Lemma~\ref{lem:formula:surprise} we have
\[D_{\text{surpr}}^\Ac=\{d\in D \cap K_\Ac \mid K_\Ac\not\subseteq\{\Mo,\dots,d\} \}.
\]
We consider two cases. 
\begin{description}
\item{Case $D_{\text{surpr}}^\Ac=\emptyset$: \ } Proposition~\ref{prop:tau:basic} gives us that $\tau^\Ac$ evaluates to $\bot$.
The description in {\SpastL} implies that $\Ac\vdash\tau^\Ac$, and hence $\Ac$ is inconsistent. 
Thus, (a) is true for $\Ac$ because it is inconsistent, 
and (b) is also true for $\Ac$, because both conditions are fulfilled.
\item{Case $D_{\text{surpr}}^\Ac\neq\emptyset$: \ } We first note that situation (b) is impossible.
(Since $L$ is an $\Ac$-p-knowledge set, all elements of $D_{\text{surpr}}^\Ac$
are different from $\max L$.)
Proposition~\ref{prop:tau:basic} says that $\tau^\Ac$ evaluates to $\top$.
By the description in {\SpastL} we have $\Ac\sim\Ac_L$, and hence $\Ac\nvdash{\langle\cT\in K\rangle}$ for all $K\subsetneq L$.
Hence $L=K_\Ac$. With ``(ii) $\Leftrightarrow$ (iii)'' in Lemma~\ref{lem:formula:surprise} we see that formula $\tau^\Ac$ shrinks to 
\[\tau^\Ac\;=\;\bigvee_{d\in D}[d\in L-\{\max L\}]\;\equiv\;[D\cap L-\{\max L\}\neq\emptyset],\]
and $D_{\text{surpr}}^\Ac=D\cap L-\{\max L\}$, by Definition~\ref{def:surprising:run:formula} and Lemma~\ref{lem:formula:surprise}(iii). This implies (a).\qed
\end{description}

One should note how the second case in the proof shows that the self-reference in {\SpL}
collapses to nothing once it is applied with a rational teacher. 
The information that the students can use is that there exists a day in $D\cap L$ on which a surprising test is possible, 
and this is a simple truth value, with no bearing on the actual choice of the teacher. 
Moreover, let us imagine that $D\cap L-\{\max L\}\neq\emptyset$ and that the students are informed 
that it is the goal of the teacher to have a surprising test.
Assuming their teacher acts rationally (what else?), they can figure out for themselves that
a surprising test will take place. This shows that the announcement carries no further information beyond what the goal of the 
teacher is. Thus we should not be astonished that it is of no use for further derivations. 
(Recall the situation with the chess player playing White who has a one-move checkmate: The Black player can see
for himself that this is so, and no announcement is needed to inform him.) 

We note that the previous remarks explain why our result fully vindicates the completely naive view,
which is that all legal choices of the teacher excepting $\max L$ will create surprise. 
This naive view works as long as the ``self-reference'' is ignored. 
(Quine's analysis~\cite{Quine:1953} is an example, and in~\cite{Shaw:1958} the effects of 
stepwise acknowledging this self-referential feature are spelled out. Later works, beginning with~\cite{Fitch:1964}
find that full self-reference will create a contradictory situation.) 
Our new perspective, namely that of taking a rational teacher into account without saying this is guaranteed, 
shows that even with full self-reference the students can prove nothing beyond the naive view.

\subsection{Central examples}
%%%%%%%%%%%%%%%%%%%%%%%%%%%%%%%%%%%%%%%%%%%%
\begin{description}

\item{$L=R$: \ } This is the setting of the original story. 
The solution is $\Ac=\Ac_R$, $K_\Ac=R$ and $D_{\text{surpr}}^\Ac=D$. 
A rational teacher chooses an arbitrary day and thus ensures that she is successful; she will avoid the choice $\none$.

\item{$L=D$: \ } This is the ``surprising egg paradox'' from~\cite{Scriven:1951}.
The solution is $\Ac=\Ac_D$, $K_\Ac=D$ and $D_{\text{surpr}}^\Ac=D\cap \{\Mo,\dots,\Th\}=\{\Mo,\dots,\Th\}$.
A rational teacher chooses an arbitrary day excepting Friday. 

\item{$L=\{\Mo,\none\}$: \ } This is the ``birthday present paradox'',
after~\cite{Gardner:1991,Scriven:1951}.
The solution is $\Ac = \Ac_L$, $K_\Ac=\{\Mo,\none\}$ and $D_{\text{surpr}}^\Ac=\{\Mo\}$.
A rational teacher chooses Monday.
Although this choice is uniquely determined for a rational teacher, it is not possible to prove $\langle\cT=\Mo\rangle$ from $\Ac_L$,
see the proof of (a), (b), (c) of Theorem~\ref{thm:students:are:wrong}.
\end{description}

The last example is very interesting, because it condenses
the question of truth or falsity of the teacher's utterance as far as possible.
We may also look back at our Chess example from Section~\ref{sec:one:move}. The Black player might even say
``I will checkmate you in one move, but you cannot know/prove it.'' 
This statement expresses the possibility of a one-move checkmate, and at the same time that there are alternative non-winning moves.
And it is completely void of information, since the White player can see for himself what is being said!

\subsection{Final remarks on ``rational teacher''}
%%%%%%%%%%%%%%%%%%%%%%%%%%%%%%%%%%%%%%%%%%%%%%%%%%%%%%%%%%
Sorensen~\cite{sorensen1988blindspots} reports, 
citing several sources, that the story is based on real events
like the announcement of a ``surprise blackout''.  
If this is so, the story was not constructed as a logical puzzle, 
and the tricks with the inaccuracies of everyday language were not built into it on purpose.
Likewise, rational teachers are not something woven into it by design. 
Still, this concept feels like a perfect fit.
On the other hand, it should be clear that it is not really necessary for understanding
what is going on and what is going wrong in the story, see Section~\ref{sec:students}. 

Based on our results, we can give the original announcement {\Sp} 
a detailed formulation in everyday language, like in the following. 
\begin{description}
\item{{\Spast} \ } \begin{minipage}[t]{0.82\textwidth}
``(i) You are certainly aware that my objective is to give a surprising test next week.\\[2pt]
	(ii) I can give a test that is ``surprising''
	in the sense that on the morning of the day with the test you cannot ``know'', i.e., prove, that it is on this day, 
	because I will have an alternative that you cannot exclude. 
	This holds even if you use for your proof that (ii) is true!\\[2pt]
	[~(iii) \ The last sentence in (ii) is just mocking you, 
	because you can infer (ii) from (i) by yourself, so (ii) is no new information.~]\\[2pt] 
	(iv) If you assume rational behavior on my side you will see that there will be a surprising test next week.'' 
\end{minipage}
\end{description}
There is no contradiction and nothing paradoxical in {\Spast}.

%%%%%%%%%%%%%%%%%%%%%%%%%%%%%%%%%%%%%%%%%%%%%%%%%%%%%%%%%%%%%%%%%%%%%%%%%%%%%
\section{The ``paradox'' and modal logic}\label{sec:modal:logic}
%%%%%%%%%%%%%%%%%%%%%%%%%%%%%%%%%%%%%%%%%%%%%%%%%%%%%%%%%%%%%%%%%%%%%%%%%%%%%
We give here a description of our argument against the ``paradoxness'' of the surprise test story once again, 
now in the language of (propositional) modal logic. The ideas are the same, but there are differences in the technical execution. 
We include this part mainly to help in understanding the relationship from 
our point of view to traditional ways of looking at STP.
Moreover, at least in the standard setting, i.e., in system S5, the deliberations become
even simpler than in plain propositional logic, since the concept of ``knowing something''
can be expressed by a modal formula in the logic itself, without tricks like Iverson brackets to ``import'' 
the truth value of provability statements. It will turn out that for the formula $\sigma$ that expresses
that there is a surprise test the argument of the students corresponds to taking $\nec\sigma$ 
as the starting point, but that this formula is unsatisfiable. 
The arguably reasonable way to understand the announcement is to take $\poss\!\sigma$ as 
the assumption. This will make it possible to accommodate ``self-reference'' and create no contradiction at all. 

\subsection{Elements of propositional modal logic}\label{sec:Kripke:basic}
%%%%%%%%%%%%%%%%%%%%%%%%%%%%%%%%%%%%%%%%%%%%%%%%%%%%%%%%%%%%%%%%%%%%%%%%%%%
Propositional modal logic extends usual propositional logic. 
For an introduction see~\cite[Ch.~2, 4, 5]{Boolos:1994} or~\cite{2008:Priest};
for more detailed information see~\cite{2007:BlackburnVanBenthem,2007:GorankoOtto}.
We compactly list the basic notions we need here, specialized to our situation.
The set of formulas of propositional modal logic (``formulas'' in this section) is defined inductively, as follows.
(i) $\bot$ is a formula. (ii) $X_i$ is a formula, for $1\le i\le 5$.
(iii) If $\varphi$ and $\psi$ are formulas then $(\varphi\to\psi)$ is also a formula. 
(iv) If $\varphi$ is a formula, then $\nec\varphi$ is also a formula.
Formula $\nec\varphi$ is read as ``box $\varphi$''. Other ways of reading it are ``$\varphi$ is necessary'', 
``it is known that $\varphi$ holds'', ``it is believed that $\varphi$'', and so on, depending on the context and
the way the logic is interpreted, see, e.g.,~\cite{2008:Priest}.%
\footnote{The technical part of our discussion does not use a specific interpretation of $\nec$.
But of course our results are open to such interpretations.}
The operators $\neg$, $\vee$, $\wedge$, $\leftrightarrow$, and $\top$ are then defined as in plain propositional logic, see Section~\ref{sec:propositional:formulas}.

Further we let $\poss\!\varphi:=\neg\!\nec\!\neg A$, read as ``diamond $A$'' or ``$A$ is possible''. 

We next introduce Kripke models (\cite{Kripke:1963}, also see~\cite{2007:GorankoOtto}), the standard tool for giving semantics to modal formulas,
with some simplifications that are admissible because there are only finitely many variables.

A \emph{Kripke model}
$\Mc=(W,E,V)$ consists of a nonempty set $W$, a binary relation $E$ on $W$
and a \emph{valuation} $V\subseteq W\times\{X_1,\dots,X_5\}$.
The elements $w$ of $W$ are called \emph{worlds}. 
If $(w,w')\in E$, we often say that $w'$ is a \emph{successor} of $w$, that world $w$ \emph{sees} world $w'$, or
that $w'$ \emph{is visible} from $w$.
It is often useful to visualize $E$ as a set of edges (``\emph{arcs}'') in a directed graph $(W,E)$.
A valuation $V$ induces, for each world $w$, 
an assignment $a_w$ for the variables, by letting $a_w(X_i):=[(w,X_i)\in V]$.% 
\footnote{See Section~\ref{sec:propositional:formulas} for an explanation of Iverson brackets $[\,\cdot\,]$.}
Of course, all assignments $a_w$ together in turn determine $V$. We will often work with the notation $a_w$. 

Truth values for formulas $\varphi$ in worlds $w$ of arbitrary Kripke models $\Mc=(W,E,V)$ are defined inductively, as follows.
\begin{align*}
&\llbracket\bot\rrbracket_{\Mc{,}w}&=&\;\;0;\\
&\llbracket X_i\rrbracket_{\Mc{,}w}&=&\;\;a_w(X_i);\\
&\llbracket(\varphi\to \psi)\rrbracket_{\Mc{,}w}&=&\;\; \max\{(1-\llbracket\varphi\rrbracket_{\Mc{,}w}),\llbracket\psi\rrbracket_{\Mc{,}w}\};\\
&\llbracket\nec\varphi\rrbracket_{\Mc{,}w}&=&\;\;\min\{\llbracket \varphi\rrbracket_{\Mc{,}w'}\mid w'\in W,(w,w')\in E\}.\hspace*{2cm} 
\end{align*}
If $\llbracket\varphi\rrbracket_{\Mc{,}w}=1$ we say ``$\varphi$ \emph{holds} in world $w$ (of $\Mc$)'' or ``$\varphi$ is satisfied in world $w$ (of $\Mc$)'', 
and write $\Mc{,}w\models \varphi$.
If $\varphi$ holds in all worlds of $\Mc$, we say ``$\varphi$ holds in $\Mc$'' and write $\Mc\models \varphi$.
\footnote{When talking about formulas and proofs, we will never say a formula is $\varphi$ ``true''. 
The problem with such formulations is that it is not clear whether $\varphi$ holds in a world, in a model, or in all worlds.} 
If $\varphi$ holds in all worlds of all models, we say that $\varphi$ is \emph{valid}, and we write $\models\varphi$.%

\begin{fact}\label{fact:gueltigkeit}
Let $\Mc$ be a Kripke model, and let $\varphi$ and $\psi$ be formulas. Then we have:  
\begin{itemize}
\item[(a)] If $\varphi$ is a propositional tautology, 
i.e., if considering partial formulas $\nec \psi$ of $\varphi$ as atomic gives a tautology, 
then $\varphi$ is valid.
\item[(b)] If $\varphi$ and $\varphi\to \psi$ both hold in world $w$, then $\psi$ holds in world $w$.
\item[(c)] $\nec(\varphi\to \psi)\to(\nec \varphi\to \nec \psi)$ is valid. 
\item[(d)] If $\Mc\models \varphi$ then $\Mc\models\nec\varphi$.\qed 
\end{itemize} 
\end{fact} 
For the simple proofs of facts like this the reader is referred to the mentioned textbooks.
The following is easily read off from the definitions. 

\begin{fact}\label{fact:and:mod:ponens:under:box}
In all Kripke models $\Mc$ we have, for all worlds $w$ and all formulas $\varphi$ and $\psi$:
\begin{itemize}
	\item[(a)] $\Mc{,}w\models \; (\nec \varphi\wedge \nec \psi) \; \leftrightarrow \; \nec(\varphi\wedge \psi)$.
	\item[(b)] If $\Mc{,}w\models \nec \varphi$ and $\Mc\models \varphi\to \psi$, then $\Mc{,}w\models \nec \psi$.
	\item[(c)] If $\Mc{,}w\models \nec \varphi$ and $\Mc{,}w\models \nec (\varphi\to \psi)$, then $\Mc{,}w\models \nec \psi$.\qed
\end{itemize}
\end{fact}

Whether $\Mc{,}w\models\varphi$ holds or not is decided within a \emph{submodel} $\Mc_w:=(W_w,E_w,V_w)$, which
consists of those worlds that can be reached from $w$ along $E$-paths, i.e.,  $W_w=\{u\in W\mid (w,u)\in E^i\text{ for some $i\ge0$}\}$, 
$E_w= E\cap W_w^2$, and $V_w=V\cap (W_w\times\{X_1,\dots,X_5\})$. 

\begin{fact}\label{fact:submodel:at:w}
For all Kripke models $\Mc$, all worlds $w$ in $\Mc$, and all formulas $\varphi$ we have
$\Mc{,}w\models\varphi\text{ \ if and only if \ }\Mc_w,w\models\varphi$.
\qed
\end{fact}
Let $\Ac$ be a set of formulas.
We write $\Mc\models \Ac$ and say ``$\Mc$ satisfies $\Ac$'' or ``$\Ac$ is satisfied in $\Mc$'',
if all formulas $\varphi \in \Ac$ hold in $\Mc$. We want to define when 
a formula $\varphi$ is a (semantic) \emph{consequence of} $\Ac$ (or \emph{follows from} $\Ac$).
There are two versions of this, one \emph{global} (in units of models), one \emph{local} (in units of worlds), of which we choose the latter one.
We define:
\[\Ac\models \varphi :\Leftrightarrow \Mc{,}w\models \varphi \text{ for all worlds $w$ in a Kripke model $\Mc$ with $\Mc{,}w\models\Ac$.}\]

The following specializations of Kripke models will be relevant for our work.
\begin{definition}
\label{def:transitiv:reflexiv}
A Kripke model $\Mc=(W,R,V)$ is called
\begin{itemize}
	\item[(a)] \emph{transitive} if $E$ is transitive, and 
	\item[(b)] \emph{reflexive} if $E$ is reflexive, i.e., if $(w,w)\in E$ for all $w\in W$, and
  \item[(c)] an \emph{equivalence model} if $E$ is transitive, reflexive, and symmetric, i.e., if $E$
is an equivalence relation. 
\end{itemize}	
\end{definition}
The following are simple facts for these specialized types of Kripke models, which 
easily follow from the definitions, see, e.g.,~\cite[Ch.~4]{Boolos:1994}.
\begin{fact}\label{fact:transitive:reflexive}
Let $\Mc$ be a Kripke model. Then we have, for all formulas $\varphi$:
\begin{itemize}
\item[(a)] $\Mc$ is transitive  \ $\Rightarrow$ \ $\Mc\models \nec\varphi\to \nec\!\nec\!\varphi$, 
 \ ( $\Leftrightarrow$ \ $\Mc\models \,\poss\!\!\poss\!\!\varphi\to \poss\!\varphi$ \ ).
\item[(b)] $\Mc$ is reflexive \ $\Rightarrow$ \ $\Mc\models \nec\varphi\to \varphi$ 
\ ( $\Leftrightarrow$ \ $\Mc\models \varphi\to \poss\!\varphi$).
\item[(c)] $\Mc$ is an equivalence model \ $\Rightarrow$ \ $\Mc\models \nec\varphi\to \varphi$
 and  $\Mc\models \poss\!\varphi\to \nec\!\poss\!\!\varphi$.

\end{itemize}
\end{fact}
We give names to these sets of formulas valid in transitive and reflexive models, and equivalence models, respectively.  
\begin{align*}
\mathcal{TR}&:=\{\nec \varphi\to \varphi\mid\text{$\varphi$ is formula}\}\cup\{\nec \varphi\to \nec\!\nec\!\varphi\mid\text{$\varphi$ is formula}\}.\\
\mathcal{E\kern-0.5ptQ}&:=\{\nec \varphi\to \varphi\mid\text{$\varphi$ is formula}\}\cup\{\poss\!\varphi\to \nec\!\poss\!\!\varphi\mid\text{$\varphi$ is formula}\}.
\end{align*}

\subsection{Systems of modal logic}\label{sec:modal:logic:proofs}
%%%%%%%%%%%%%%%%%%%%%%%%%%%%%%%%%%%%%%%%%%%%%%%%%%%%%%%%%%%%%%%%%%%%%%%%%%%%%%%%%%%%%%
There are several proof calculi for modal logic, each given by an axiom system.
Our basis is the simple variant that is called 
\emph{system K}, which has the following axioms:
\begin{itemize}
	\item[$(\alpha)_{\text{K}}$] all formulas $\varphi$ that are propositional tautologies, 
	\item[$(\beta)_{\text{K}}$] $\nec(\varphi\to\psi)\to (\nec\varphi \to \nec\psi)$, for formulas $\varphi$ and $\psi$. 
\end{itemize}

By Fact~\ref{fact:gueltigkeit}(a) and (c) these axioms hold in all worlds in all Kripke models. 
The set of formulas provable in system K is defined inductively, as follows.                                          
\begin{itemize}
\item[(i)] $\varphi$ is provable for all axioms from $(\alpha)_{\text{K}}$ and $(\beta)_{\text{K}}$.
	\item[(ii)] \emph{modus ponens}: If $\varphi$ and $\varphi\to\psi$ are provable then $\psi$ is also provable.
	\item[(iii)] \emph{necessitation}: If $\varphi$ is provable then $\nec\varphi$ is also provable.
\end{itemize}

We write \mbox{$\vdash\varphi$} if $\varphi$ is provable in system K. 
By Fact~\ref{fact:gueltigkeit}(a), (b), and (d), basis (i) and the \emph{deduction rules} (ii) and (iii) are \emph{correct} for Kripke models,
in that for each Kripke model $\Mc$ from formulas that hold in $\Mc$ these rules always yield formulas that also hold in $\Mc$. 
Clearly, then, all formulas $\varphi$ with \mbox{$\vdash\varphi$} hold in all Kripke models. 

Now consider an arbitrary set $\Ac$ of formulas. We define the set of formulas \emph{K-provable from} $\Ac$ 
as the smallest set of formulas that contains $\Ac$ and is closed under (i)--(iii).
We write $\Ac\vdash \varphi$ if $\varphi$ is K-provable from $\Ac$.
Trivially, $\emptyset\vdash \varphi$ means the same as $\vdash \varphi$.
The following fact summarizes the connection between K-provability and
semantic consequence in arbitrary Kripke models. For a proof see, e.g.,~\cite{Boolos:1994,2008:Priest}.
 
\begin{fact}[\mbox{(Correctness and completeness of K-provability)}]\label{fact:K:correctness}
For each set $\Ac$ of formulas and every formula $\varphi$ we have $\Ac \vdash \varphi$ $\Leftrightarrow$ $\Ac\models \varphi$.\qed 
\end{fact}

An \emph{axiom system} $\Ac$ is just a set of formulas. However, in connection with axiom systems
we are interested in the stronger version of semantic consequence, which works in units of models.
We will need this only for the axiom systems that we get if we add 
 $\mathcal{TR}$ and $\mathcal{E\kern-0.5ptQ}$, resp., from Section~\ref{sec:Kripke:basic} to K. The systems we get are called S4 and S5, resp.
System S5 is at least as strong as S4 (even strictly stronger, but we do not need this), because $\mathcal{E\kern-0.5ptQ}\vdash \varphi$ for all $\varphi$ in $\mathcal{TR}$~(see, e.g.,~\cite{Boolos:1994}). 
Accordingly, we get stronger provability notions, for which we write $\svdash$ resp. $\sfdash$. Just as in Fact~\ref{fact:K:correctness}
we have~\cite{Boolos:1994,2008:Priest}:

\begin{fact}[(Correctness and completeness for S4 and S5)]\label{fact:korrektheit:transitive:reflexiv}
For every formula $\varphi$ we have:
\begin{description}
	\item{$\svdash \varphi$} \ $\Leftrightarrow$ \ $\Mc\models \varphi$ for all transitive and reflexive models $\Mc$ of $\Ac$.
	\item{$\sfdash \varphi$} \ $\Leftrightarrow$ \ $\Mc\models \varphi$ for all equivalence models $\Mc$ of $\Ac$.\qed 
\end{description}
\end{fact}

Finally, we define a set $\Ac$ of formulas to be called \emph{consistent} if there is a formula $\varphi$ with $\Ac\not\vdash \varphi$; 
otherwise $\Ac$ is \emph{inconsistent}. Obviously, $\Ac$ is consistent if and only if $\Ac\not\vdash \bot$.
Fact~\ref{fact:K:correctness} gives that $\Ac$ is consistent if and only if there is a world $w$ in a Kripke model $\Mc$ with $\Mc{,}w\models\Ac$.
In the same way, this fact implies that $\Ac\not\vdash \varphi$ if and only if there is a Kripke model $\Mc$ with a world $w$ for which $\Mc\models\Ac$ and $\Mc{,}w\models \neg\varphi$.  
The situation is completely analogous for S4- and S5-provability from formula set $\Ac$ 
and the corresponding narrower model classes, with Fact~\ref{fact:korrektheit:transitive:reflexiv}.

\begin{remark}\label{bem:induced:submodel:sfour:sfive}
Let $\Mc$ be a transitive and reflexive model, let $w$ be a world in $\Mc$, and let $\varphi$ be a formula.
We have noted after Fact~\ref{fact:and:mod:ponens:under:box} that the truth value $\llbracket\varphi\rrbracket_{\Mc{,}w}$ only depends on
the submodel $\Mc_w:=(W_w,E_w,V_w)$, which consists of the worlds reachable from $w$ along $E$-paths.
This is particularly simple in reflexive and transitive models, where $W_w$ is the set of worlds visible from $w$.
(If $\Mc=(W,R,V)$ is an equivalence model, $W_w$ is the equivalence class of $w$.)
In this case we have: 
\begin{equation}\label{eq:von:box:zu:forall}
\Mc,\!w\models \nec\varphi  \;\Leftrightarrow\; \forall u{\in}\, W_w\colon\Mc_w,\!u\models\varphi\;\Leftrightarrow\; \Mc_w\models\varphi.
\end{equation}
Negating the formulas one gets:
\begin{equation}\label{eq:von:raute:zu:exists}
\Mc,\!w\models \poss\!\varphi \;\Leftrightarrow\; \Mc_w,\!w\models \poss\!\varphi \;\Leftrightarrow\;  \Mc_w\not\models \neg\varphi.
\end{equation}
\end{remark}

\subsection{Kripke models for STP}
%%%%%%%%%%%%%%%%%%%%%%%%%%%%%%%%%%%%%%%%%%%%%%%%%%%%%%%%
Just as in the case of propositional logic (see Section~\ref{sec:propositional:formulas}), 
we adopt axiom (Ax$_{\le1}$), define auxiliary variables
$Y_\Mo,\dots,Y_\Fr$, and $Y_\none$ to be used in the following, and identify
(Ax$_{=1}$) as a consequence of (Ax$_{\le1}$). This is tacitly assumed as axiom for everything we do from here on.
Other notation defined in Section~\ref{sec:propositional:formulas} (e.g., formulas like  $\langle\cT=d\rangle$) is used without further comment.   

There are many different Kripke models for STP.
The model in which for each run $r$ 
there is exactly one $r$-world $w_r$ 
and in which all worlds ``see'' each other (including itself) is of particular interest.
An important variant is a model in which run $\none$ does not occur. 
\begin{definition}\label{def:the:model}
\begin{itemize}
	\item[(a)]
	      $\McG=(\WG,\EG,\VG)$ with $\WG=\{w_r\mid r\in R\}$, $\EG=\WG\times \WG$, and
        $V_G$ so that $a_{w_r}(Y_r)=1$ for $r\in R$.	
	\item[(b)] 
	      $\McG^+=(\WG^+,\EG^+,\VG)$ with $\WG^+=\{w_d\mid d\in D\}$, $\EG^+=\WG^+\times\WG^+$, and
        $V_G$ so that $a_{w_d}(Y_d)=1$ for $d\in D$.
\end{itemize}
\end{definition}
Obviously, $\McG$ and $\McG^+$ are equivalence models, in which the visibility relations $\EG$ and $\EG^+$, resp., have exactly one equivalence class. 
We obtain similar models from $\McG$ and $\McG^+$ by leaving out one or more $d$-worlds, for $d\in D$.
Altogether there are 63 models of this kind (because the set of worlds cannot be empty).

The general structure of a model $\Mc=(W,E,V)$ for STP is the following:
$W$ is a nonempty set, each world $w\in W$ is an $r$-world for some run $r$. 
Relation $E$ is transitive%
\footnote{We choose the restriction of transitivity to model the following natural idea. 
If $(w,w')\in E$ and the students ``sit'' in world $w$, they ``see'' world $w'$
including the full situation in $w'$, i.e., all formulas that hold or do not hold in $w'$.
Now if an $s$-world $w''$ is visible from $w'$ then $\Mc{,}w'\models \poss Y_s$. 
We would like to have that now not only \mbox{$\poss\kern-1.4pt \poss\kern-1.4pt Y_s$} holds in $w$, but $\poss Y_s$ directly.
Exactly this is achieved by the requirement that $E$ be transitive.}
and reflexive%
\footnote{Reflexivity models the idea that the students perceive the world in which they ``sit'' as a possible world.}%
.
If non-symmetric visibility relations $E$ are admitted, 
one should not think of these models as too simple. 
There could be several worlds for one run $r$, with different properties, and relation $E$ is an arbitrary quasiorder.
Here and in equivalence models, a model can have several components that are not connected by $E$-edges.

In contrast, for equivalence models and S5, the situation is rather simple.%
\footnote{This is so not only for the STP situation, but for 
systems with finitely many propositional variables in general.}
For an equivalence model $\Mc$ and a world $w$ in $\Mc$ let $B_w:=\{s\in R\mid\text{$W_w$ has an $s$-world}\}$.
Now consider some model $\Mc'$ and a world $w'$ such that $B_w=B_{w'}$.
One can then prove by induction on formulas that for all $r\in B_w$
and all $r$-worlds $u$ in $W_w$ and all $r$-worlds $u'$ in $W_{w'}$
a formula $\varphi$ holds in $u$ if and only if it holds in $u'$.
This entails that the set of formulas that hold in a world $w$ in an equivalence model
are determined by $B_w$ and $r_w$, where $w$ is an $r_w$-world. 
From this we can conclude that for each pair $(B,r)$ with $r\in B\subseteq R$
all worlds $w$ with $(B_w,r_w)=(B,w)$ satisfy exactly the same formulas. 
This makes it possible to build an explicit \emph{universal model}. 
\begin{definition}\label{def:universal:model}
 $\Mc^*:=(W^*,E^*,V^*)$, where
\begin{align*}
W^*&=\{(B,r)\mid r\in B\subseteq R\},\\
E^*&=\{(B,r)(B',r')\in W^*\times W^*\mid B=B'\},\\
V^*&=\{((B,d),X_d)\mid (B,d)\in W^*, d\in D\}.
\end{align*}
\end{definition}
This definition entails that $(B,r)$ is an $r$-world,
and that the equivalence class $W_w$ of $w=(B,r)$ has an $s$-world if and only if $s\in B$.

With Fact~\ref{fact:korrektheit:transitive:reflexiv}
it is routine to show the following. 
\begin{observation}\label{obs:univeral:model}
$\sfdash\varphi$ if and only if $\Mc^*\models\varphi$, for all formulas $\varphi$.
\end{observation}

Clearly, the universal model is finite (it has $6\times 2^5=192$ worlds and 63 equivalence classes). 
Quite trivially, then, just checking $\Mc^*,w\models\varphi$ for all worlds $w\in W^*$
gives an algorithm to decide whether $\varphi$ is S5-valid. 
Moreover, it is possible to prove that all that can be expressed by an arbitrary formula set
can be condensed in one single formula from a suitable finite set of formulas. 

\begin{lemma}\label{lemma:finiteness:in:equivalence:models}
There is a finite set $P$ of formulas such that for 
each formula set $\Ac$ there is some $\tau_\Ac$ in $P$ with 
$\Ac\sfdash\tau_\Ac$ and $\sfdash\tau_\Ac\to\varphi$ for all $\varphi$ in $\Ac$.
\end{lemma}
\emph{Proof}:
Consider \[W^*_\Ac=\{w\in W^*\mid\Mc^*{,}w\models\varphi\text{ for all $\varphi$ in $\Ac$}\}.\]
In view of what was just said about $\Mc^*$ we need a formula $\tau=\tau_\Ac$ with $W^*_\Ac=W^*_{\{\tau\}}$. The following definition is suitable.
\[\tau_\Ac:=\bigvee_{(B,r)\in W^*_\Ac}\!\!\langle\cT=r\rangle\wedge\nec\langle\cT\in B\rangle\wedge\bigwedge_{s\in B}\!\poss\!\langle\cT=s\rangle.\]
As there are only finitely many subsets of $W^*$, there are only finitely many different formulas $\tau_\Ac$.
\qed

\subsection{``Surprise'' in modal logic}\label{sec:surprise:modal}\label{sec:surprise:in:modal:logic}
%%%%%%%%%%%%%%%%%%%%%%%%%%%%%%%%%%%%%%%%%%%%%%%%%%%%%%%%%%%%%%%%%%%%%%%%%%%%%%%%%
Modal logic makes it possible to express provability (or ``knowledge'', or ``belief'') directly as a formula. 
In contrast to the situation in propositional logic we do not need tricks like Iverson brackets to force such statements into the world of formulas.

We first consider the question of ``surprise'' relative to a fixed world $w$ in a model $\Mc$, 
and with the choice of the teacher being $\rT\in R$.
We may assume that $w$ is an $\rT$-world, since otherwise it is irrelevant what is the matter in $w$. 
If $\rT=\none$, there is no test in $w$, so certainly not a surprising test.
So assume $\rT=d\in D$. This means $\Mc{,}w\models\langle \cT= d\rangle$.
If $\Mc{,}w\models\nec \langle \cT\le d\rangle$ (from $w$ only $r$-worlds with $r\le d$ are visible)
then in world $w$ the students ``know'' or ``believe'' that no run $r>d$ is possible, and we should certainly say that
in $w$ they ``expect'' that there is a test on day $d$. What about the reverse direction?

Let us imagine the situation on the morning of day $d$. The students may add an axiom $\langle \cT\ge d\rangle$
to what they see in $w$, which effectively removes all $d'$-worlds with $d'<d$ from $W_w$, leaving a submodel $\Mc_w|_{\ge d}$.
Can we now say that if some $r$-world $w'$ with $r>d$ exists in $\Mc_w|_{\ge d}$ then the test on day $d$ in $w$ is ``surprising''?
Curiously, this is not necessarily so in general transitive and reflexive models
(corresponding to S4). Namely, it could be the case that $(\Mc_w|_{\ge d},w)\nesi(\Mc_{w'}|_{\ge d},w')$ does not hold, 
and in this case the students will ``know'' that they are not in the $r$-world $w'$.%
\footnote{It is not very hard to construct an example for such a situation, but it is not of too much interest.} 
So to define ``surprise'' before the background of Kripke models (for S4) we need to be more careful.
For the $d$-world $w$ in $\Mc$ a test on day $d$ should be ``surprising'' if and only if
the students ``cannot know'' on the basis of $\Vn(\Mc|_{\ge d},w)$ 
whether they are in world $w$ with a test on day $d$ or in some ``${>}d$-world'', i.e., an $r$-world $w'$ with $r>d$.

\begin{definition}\label{def:surprise:day:d:in:w}
  Let $\Ac$ be a set of formulas and let $w$ be a world in an equivalence model $\Mc$. 
	We say that in world $w$ there is an \emph{$\Ac$-surprising test on day $d$} with respect to S5 
	if $w$ is a $d$-world with $\Mc{,}w\models\Ac$, 
	and there is an equivalence model $\Mc'$ with a ${>}d$-world $w'$ such that $(\Mc|_{\ge d},w)\nesi(\Mc'|_{\ge d},w')$.\\[2pt]
(The same definition is applied for S4, replacing ``equivalence'' by ``transitive and reflexive''.)
\end{definition}

Beware that it is \emph{not} required that $\Ac$ is also satisfied in $w'$.
Already this choice reflects much of the deliberations of Section~\ref{sec:one:move}
where it is explained that we should assume rational behavior on the side of the teacher. 
If she has a chance to make good on her announcement she will do it even if 
other options are available to her that do not fulfill the announcement. 
 
By the restrictions $\Mc|_{\ge d}$ and $\Mc'|_{\ge d}$ it is captured that the days before 
$d$ have passed, and by the $\nec$-equivalence that it is impossible to see a difference between $w$ and $w'$
by looking at arbitrary $\nec$-formulas that hold in these two worlds.

\begin{example}\label{example:standard:model}
In $\McG$ all worlds in $\WG$ ``see'' all worlds in $\WG$, 
and hence they all have the same set $\Vc_{\nec}(\Mc|_{\ge d},w)$,
which expresses that all worlds $w_r$ with $r\ge d$ are visible, 
in particular the $\none$-world. In formulas: $\Vc_{\nec}(\McG|_{\ge d},w)$ is the set of all S5-consequences of
\[\nec\langle\cT\ge d\rangle\;\wedge\;\bigwedge_{r\ge d}\poss\!\langle\cT=r\rangle,\]
for all worlds $w$ in $\McG|_{\ge d}$. The same is true for $w_\none$, and we have $\none>d$. 
This shows that for $d=\Mo,\dots,\Fr$ in $w_d$ there is 
a test on day $d$ that is $\emptyset$-surprising (or \emph{surprising}).
\end{example}

With Definition~\ref{def:surprise:day:d:in:w} we have described \emph{mathematically} when there is
an $\Ac$-surprising test on day $d$ in world $w$ of a model $\Mc$ (equivalence or transitive and reflexive). 
The obvious disadvantage of this definition is that it is complicated and that it uses the whole
semantics machinery of modular logic. There is no obvious way how this condition can be checked algorithmically.
Thus we would like to find a formula $\sigma^\Ac_d$ that characterizes worlds with this property. 
We shall presently see that this is quite easy in equivalence models. 
For general transitive and reflexive models a somewhat larger technical effort is necessary.  

%%%%%%%%%%%%%%%%%%%%%%%%%%%%%%%%%%%%%%%%%%%%%%%%%%%%%%%%%%%%%%%%%%%%%%%%%%%%%%%%%%%%%%%%%%%%%%%%%%%%
\subsection{``Surprising test'' in equivalence models}\label{sec:surprise:sf}
%%%%%%%%%%%%%%%%%%%%%%%%%%%%%%%%%%%%%%%%%%%%%%%%%%%%%%%%%%%%%%%%%%%%%%%%%%%%%%%%%%%%%%%%%%%%%%%%%%%%
In this section, by ``model'' we always mean an equivalence model,
and the proof calculus uses S5. 
From Lemma~\ref{lemma:finiteness:in:equivalence:models} we have that in the situation of Definition~\ref{def:surprise:day:d:in:w}
it is sufficient to consider formula sets $\Ac=\{\tau\}$. 
(We write ``$\tau$-surprising'' insted of the pedantic ``$\{\tau\}$-surprising''.)

\subsubsection{``Surprise'' in single worlds}
%%%%%%%%%%%%%%%%%%%%%%%%%%%%%%%%%%%%%%%%%%%%%%%%%%%%%%%%%%%%%%%%
\begin{theorem}[(``$\tau$-surprising test on day $d$'' as a formula)]\label{thm:surprise:eq:model}
Let $\Mc=(W,E,V)$ be an equivalence model, let $w\in W$, and let $\tau$ be a formula.
For each day $d$ the following are equivalent: 
\begin{itemize}
	\item[(i)] In $w$ there is a $\tau$-surprising test on day $d$, in the sense of Definition~\ref{def:surprise:day:d:in:w}. 
	\item[(ii)] $\Mc{,}w\models \tau$ and $\Mc|_{\ge d}{,}w\models\langle\cT=d\rangle  \wedge \neg\nec\!\langle\cT=d\rangle$.
	\item[(iii)] $\Mc{,}w\models  \tau\wedge\langle\cT=d\rangle \wedge \neg\nec\!\langle\cT\le d\rangle \;\;\bigl(\;\equiv\;
	 \tau\wedge\langle\cT=d\rangle  \wedge\; \poss\!\langle\cT>d\rangle\;\bigr)$.
\end{itemize}
\end{theorem}

\emph{Proof}: ``(i) $\Rightarrow$ (ii)'': Assume (i). Then $w$ is a $d$-world, in which $\tau$ is satisfied, and there is a model $\Mc'$ of $\Ac$ with a ${>}d$-world $w'$
so that $(\Mc|_{\ge d},w)\nesi(\Mc'|_{\ge d},w')$.
Since $\Mc'|_{\ge d}$ is reflexive, we have $\Mc'|_{\ge d},w'\models\poss \neg \langle\cT=d\rangle$. 
Because $(\Mc|_{\ge d},w)\nesi(\Mc'|_{\ge d},w')$,
formula $\poss \neg \langle\cT=d\rangle $ holds in $\Mc|_{\ge d},w$ as well, or, equivalently, 
$\Mc|_{\ge d}{,}w\models\neg\nec\!\langle\cT=d\rangle$.
This establishes (ii).%
\footnote{One can go one step further: If this ${>}d$-world is $w''$, we have $(w,w''),(w'',w)\in E$.
Hence $(\Mc|_{\ge d},w)\nesi(\Mc|_{\ge d},w'')$, and that means that the conditions in Definition~\ref{def:surprise:day:d:in:w}
are realized in two worlds $w$ and $w''$ in the same equivalence class of $W$.} 

 ``(ii) $\Rightarrow$ (iii)'': Assume (ii). That means that from $w$ some ${>}d$-world is visible, in $\Mc|_{\ge d}$, hence 
the same is the case in $\Mc$.
This implies (iii).

 ``(iii) $\Rightarrow$ (i)'': Assume 
$\Mc, w\models \tau\wedge \langle\cT=d\rangle \wedge\poss\!
\langle\cT>d\rangle$.
It follows that in $\Mc$ there is a ${>}d$-world $w'$ with $(w,w')\in E$. 
Since $\Mc$ is an equivalence model, the set of worlds visible from $w$ 
and the set of worlds visible from $w'$ are the same. 
This implies $(\Mc|_{\ge d},w)\nesi (\Mc|_{\ge d},w')$, and hence (i).
\qed 

Note that (iii) only refers to $\Mc$, without the restriction ``${\ge}d$\kern1pt''.
Thus we see again that time does not play a role in STP, cf. the ``designated student paradox'' in~\cite{sorensen1988blindspots}.

The theorem says that with \[\sigma_d\;:=\;\langle\cT=d\rangle \wedge \neg\nec\!\langle\cT\le d\rangle \;\;\Bigl(\equiv\;\langle\cT=d\rangle \wedge \poss\!\langle\cT> d\rangle\Bigr)\]
we have a formula that expresses that there is a surprising test on day $d$.
Thus, formula
\begin{equation}\label{eq:sigma:plus}
\sigma\;:=\; \langle\cT\le\Fr\rangle \wedge\bigwedge_{d\in D}\!(\langle\cT=d\rangle \to \sigma_d) \;\equiv\; \bigvee_{d\in D}\! \sigma_d
\end{equation}
satisfies that for any world $w$ in any model $\Mc$ we have:%
\footnote{Gerbrandy~\cite{Gerbrandy:2007} develops the ``$\bigvee$-form'' of $\sigma\wedge\nec\bot$ for ``there is a surprising test'', 
but without reference to S5. Beware that if symmetry in $E$ is not guaranteed, $\sigma$ need not correspond to ``surprise'' in the sense of Definition~\ref{def:surprise:day:d:in:w}.}
 \begin{center}
 in $w$ there is a $\tau$-surprising test on some day \ $\Leftrightarrow$ \ $\Mc{,}w\models \tau\wedge \sigma$.
\end{center}
Now let us assume the teacher makes an announcement, which can be expressed by some formula $\alpha=\tau\wedge\sigma$.
This is what she wants to achieve and this is what she announces: $\alpha$ will happen, ``surprise'' included.

\subsection{The argument of the students}\label{sec:modal:students}
%%%%%%%%%%%%%%%%%%%%%%%%%%%%%%%%%%%%%%%%%%%%%%%%%%%%%%%%%%%%%%%%%%%%
We consider here the standard situation from the story. 

The teacher announces that there will be a surprising test, 
which we take as saying that $\alpha=\sigma\equiv\bigvee_{d\in D}\sigma_d$ will hold for her choice $\rT$.
As we have used system S5 for justifying that $\sigma$ is the right formalization of the concept of 
a surprising test, it does not make sense to consider weaker logics than S5 here.

How should the students use this announcement? 
It is difficult to utilize $\sigma$ without knowing in which equivalence class it is to be applied, 
and which subformula $\sigma_d$ will be satisfied. 
One possibility, chosen by many authors, is to consider $\sigma$ an axiom 
and restrict attention to models in which $\sigma$ holds in all worlds.
By~\eqref{eq:von:box:zu:forall}, this has the same effect as to say that $\nec\sigma$ should hold
(in the world in which the formulas are evaluated, whatever it is).

Unfortunately, as was observed many times (see, e.g.,~\cite{1968:Binkley:Surprise:Examination}) such a choice
leads into a contradiction, because $\nec\sigma$ can hold in no world. 

\begin{theorem}\label{thm:students:contradiction}
\begin{itemize}\item[]$\svdash\neg\nec\!\sigma$.\end{itemize}
\end{theorem}

\emph{Proof}:
We follow the argument of the students. The proof could be carried out by applying
the deduction rules of S4 (as has been demonstrated, e.g., in~\cite{1968:Binkley:Surprise:Examination}), 
but for a change we choose the easy way and invoke Fact~\ref{fact:korrektheit:transitive:reflexiv},
which allows us to argue via models.
We assume for a contradiction that there is a transitive and reflexive model $\Mc$ with 
a world $w$ such that $\Mc{,}w\models\nec\sigma$.
By~\eqref{eq:von:box:zu:forall}, this entails $\Mc_w\models\sigma$.
So we now have a Kripke model $\Mc=\Mc_w$ with $\Mc\models\sigma$. 
We show that this is impossible.

\textbf{Initial step:} \ Since $\langle\cT\in D\rangle$ is a factor in $\sigma$,
model $\Mc$ has no $\none$-world.

\textbf{Induction (``no-Friday'') step:} \ Fix $d\in\{\Tu,\dots,\Fr\}$, and assume that $\Mc$ has no $r$-world with $r>d$.
Since $\Mc\models\sigma$, we get $\Mc\models\bigvee_{d'\le d}\sigma_{d'}$, and hence $\Mc\models\nec\bigvee_{d'\le d}\sigma_{d'}$,
by Fact~\ref{fact:transitive:reflexive}, in particular $\Mc\models\nec\langle\cT\le d\rangle$.
Clearly, then, $\Mc\models\neg\sigma_d$, and hence $\Mc\models\bigvee_{d'\le d-1}\sigma_{d'}$.
In particular, model $\Mc$ has no $r$-world with $r>d-1$.

\textbf{The final step:} \  The result of the induction is that $\Mc$ has no $r$-world with $r>\Mo$,
so $\Mc$ can have only $\Mo$-worlds. 
On the other hand we have$ \Mc\models\sigma_\Mo$, or $\Mc\models\langle\cT=\Mo\rangle\wedge\poss\!\langle\cT>\Mo\rangle$. 
But this is impossible, which is the desired contradiction. \qed
 
We note that nowhere in Theorem~\ref{thm:students:contradiction} it is used that the announcement of the teacher is self-referential. 
This means that for the contradiction appearing this property cannot be essential.  
Statements that resemble Theorem~\ref{thm:students:contradiction} 
can be found in several works on STP (see in particular~\cite{2025:Baltag:topological:epistemic,1968:Binkley:Surprise:Examination,Scriven:1951}).
They essentially say that in the presence of $\mathcal{TR}$ ``it cannot be known that `surprise will happen'\kern1pt''.
The theorem leaves us with the question whether $\nec\sigma$ is the correct formalization of the announcement of the teacher. 
If we think this is so, we would need to step back behind S4 as our deduction system. 
But this would crush our justification for (the quite elegant and intuitive) formula $\sigma$ as the formalization of ``surprise test''. 
Thus, our plan is to salvage using S5, by arguing that it is inadequate to choose $\nec\sigma$ 
as the formalization of what the students can use.

\subsection{Capturing ``self-referentiality''}\label{sec:modal:existential:reading}
%%%%%%%%%%%%%%%%%%%%%%%%%%%%%%%%%%%%%%%%%%%%%%%%%%%%%%%%%%%%%%%%%%%%%%%%%%%%%%%%%%%%%%%%%%%%%%%%%%%%%%
Just as in Section~\ref{sec:announcement:self} we have to address the ``self-reference'' property of the announcement of the teacher.
It declares that ``surprise'' will happen, even though the students can use in their proofs that the announcement is ``true''. 
Finding an appropriate formulation for this in modal logic is a main difficulty.
As S5 is stronger than pure propositional logic, we strive hard to avoid tricks like Iverson brackets
as were used with plain propositional logic. 
Rather the announcement in its full power should be expressed as a formula 
$\alpha=\sigma\wedge\tau$ in propositional modal logic. What are the requirements for $\alpha$?
\begin{itemize}
	\item[(i)] $\alpha$ should include that there is a ``surprise test'', i.e., we require 
\end{itemize}
\begin{equation}\label{eq:alpha:implicit:1}
\alpha\;\equiv_{\text{S5}}\;\sigma\wedge\alpha\text{ \ (or \ }\sfdash\;\alpha\leftrightarrow\sigma\wedge\alpha\text{, \ or \ }\sfdash\;\alpha\to\sigma). 
\end{equation}
\begin{itemize}
	\item[(ii)] We need to capture the ``self-reference feature'' of the announcement. This means that whenever some formula $\beta$ can be shown to hold,
	mathematically, no handles barred, from the information that $\alpha$ is realized in a world $w$, this formula $\beta$ must also hold in $w$.
	By completeness of the system S5 for equivalence models this does not mean more than $\sfdash \alpha \to \text{``information that $\alpha$ holds''}$.
	But what is a formula that expresses ``the information that $\alpha$ holds''? 
	We try ``the information that $\alpha$ can happen'', or ``$\poss\!\alpha$ holds'', 
	on the basis of the idea that if $\poss\!\alpha$ holds, the teacher, being rational, will make sure that $\alpha$ will become true by her choice. 
	So our requirement becomes
	\end{itemize}
\begin{equation}\label{eq:alpha:implicit:2}
\sfdash \alpha\to\poss\!\alpha.
\end{equation}
Now~$\alpha\to\poss\!\alpha$ is an axiom of S4, hence a consequence of S5, so~\eqref{eq:alpha:implicit:2} is always true, 
and~\eqref{eq:alpha:implicit:1} is all that is required of $\alpha$.
This means that in her announcement the teacher can state an arbitrary formula $\alpha=\sigma\wedge\tau$, as long as it is not S5-contradictory,
and can make sure by her choice that $\alpha$ is realized, including ``$\alpha$-surprise'', even though $\alpha$ includes the information that $\alpha$ will happen.
(This effect is taken care of by S5, which corresponds to equivalence models, which are in particular reflexive.)
In the standard STP situation, we simply have $\alpha=\sigma$, and this is certainly not contradictory, since 
it holds in all $d$-worlds of $\McG$, for $d\in D$, see Example~\ref{example:standard:model}.

We remark a syntactic similarity of $\poss\!\alpha$to the formulas used in~\cite{2025:Baltag:topological:epistemic}, 
but note that the topological approach used there collapses in our context, since 
the topology given by equivalence model and by $\Mc^*$ is trivial in that
the worlds of an equivalence class have just this equivalence class as their only environment. 

\subsubsection{Does this work?}\label{explanation:in :model}
%%%%%%%%%%%%%%%%%%%%%%%%%%%%%%%%%%%%%%%%%%%%%%%%%%%%%%%%%%%%%%%%%
Let us look at the situation from the perspective of the univesal model $\Mc^*$,
to see what happens. Recall Observation~\ref{obs:univeral:model}. 
The teacher announces that her goal is $\alpha$, with $\sfdash\alpha\to\sigma$ and $\not\sfdash\neg\alpha$.
These properties of $\alpha$ are readily verified by the students.
As the teacher is rational, she will make a choice $\rT$ that will realize $\alpha$.
The set of possible choices is
\[D^\alpha:=\{d\in D\mid\;\not\sfdash\neg(\alpha\wedge\langle\cT=d\rangle)\},\]
and this set can also easily be calculated by the students. It is nonempty by the assumption $\not\sfdash\neg\alpha$,
and because $\sfdash\alpha\to\langle\cT\in D\rangle$.

Announcement $\alpha$ together with the idea that the teacher does not lie
lets the students (and us) focus on 
\[W^*_{\poss\!\alpha}=\{w\in W^*\mid\Mc^*{,}w\models\poss\!\alpha\},\]
i.e., worlds in those equivalence classes of $\Mc^*$ in which $\alpha$ is realized in some world.

The days $d=\Mo,\dots,\Fr$ roll by, one after the other. 
On the morning of day $d$ the students check if $d\in D^\alpha$. If this is not the case,
we have $\sfdash\neg(\alpha\wedge\langle\cT=d\rangle)$, and hence no $d$-world can satisfy $\alpha$. 
The teacher, being rational, will not choose $\rT=d$, and we will let day $d$ go. Now assume $d\in D^\alpha$.
We are especially interested in those equivalence classes of $\Mc^*$ that have a ${\ge}d$-world in which $\alpha$ holds;
so let $W^*_{\alpha,{\ge}d}$ be the set of worlds $w$ with $\Mc^*{,}w\models\poss\!(\alpha\wedge\langle\cT\ge d\rangle)$.
Assuming we are in an equivalence class inside this submodel, can the students prove that the teacher has no alternative to choosing $\rT=d$, 
assuming $\alpha$ is realized? Certainly not in equivalence classes in $W^*_{\alpha,{\ge}d}$ without a $d$-world,
and so we forget about these (if they exist), and concentrate on $W^*_{\alpha,=d}$, 
the set of worlds that satisfy $\poss\!(\alpha\wedge\langle\cT=d\rangle)$.
As $d\in D^\alpha$, we have $W^*_{\alpha,=d}\neq\emptyset$. Can the students prove that the teacher \emph{must} choose $\rT=d$,
if she \emph{does} choose $\rT=d$ and thus places the students in such a $d$-world $w$?% 
\footnote{As already mentioned, Vinogradova~\cite{2023:Vinogradova} stresses the difference between certainty and possibility, in the 
framework of constructive mathematics, where $\neg\neg\alpha$ means that $\alpha$ is possible.}
No, because $\Mc^*{,}w\models\alpha$, hence $\Mc^*{,}w\models\sigma$, hence $\Mc^*{,}w\models\sigma_d$, and hence $\Mc^*{,}w\models\poss\!\langle\cT>d\rangle$, 
which means that from $w$ an $r$-world is visible for some $r>d$, and this $r$ is a choice equally available to the teacher. 
In different equivalence classes in $W^*_{\alpha,=d}$ these alternative runs may be different, 
but every equivalence class has at least one such alternative run. 
Thus when the teacher hands out the test papers on day $d$ the students must admit they had no chance to 
prove that this must be so: No matter which equivalence class of $W^*_{\alpha,=d}$ they sit in, 
there is always an alternative choice available to the teacher.
 
In terms of provability in S5 this deliberation can be interpreted as follows.
Having been told $\alpha$, the students take $\poss\!\alpha$ as an axiom. 
They calculate $D^\alpha$, which is nonempty, and since they assume the teacher can be trusted to realize her announcement, 
they exclude all days $d\notin D^\alpha$ from consideration.
Now assume the teacher chooses $\rT=d\in D^\alpha$. When the morning of day $d$ has been reached,
the students have their axiom $\poss\!\alpha$, and they can
sharpen it to $\not\sfdash\neg\poss\!\!(\alpha\wedge\langle\cT= d\rangle)$
since the latter assertion is a consequence of $d\in D^\alpha$.
Nothing that the teacher has said or what has been observed enlarges the information the students have.
Now the test appears, and the students can see that they are in some $d$-world. 
They might not know what the equivalence class is, but this does not matter:
In each equivalence class in which $\poss\!(\alpha\wedge\langle\cT= d\rangle)$ holds, 
formula $\alpha$ holds in the $d$-world $w$. So taking the announcement into account, interpreting it as $\poss\!\alpha$, 
and fixing day $d\in D^\alpha$ as $\rT$ makes sure that $\alpha$ holds, just as the teacher has announced. In formulas:
\[\sfdash\;\;\poss\!(\alpha\wedge\langle\cT=d\rangle)\wedge\langle\cT=d\rangle \; \to \; \alpha,\text{ for $d\in D^\alpha$}.\]
``I will choose a day $\rT=d$ such that $\alpha$ will be satisfied, in particular $\sigma$ will be satisfied,
and my announcement $\alpha$ already takes all consequences of my making this announcement into account.''

\subsubsection{Examples}\label{examples:A:E}
%%%%%%%%%%%%%%%%%%%%%%%%%%%%%%%%%%%%%%%%%%%%%%%%%%%%%%%%%%%%%%%%%%%%%%%%%%

(A) $\alpha=\top$. -- The teacher does not say anything. 
This formula does not satisfiy our requirements for $\alpha$, since $\top\to\sigma$ does not hold.
Let us assume that the students have the global information
that the teacher has the goal of having a surprising test, and consider $\alpha=\sigma$. 
They can calculate that the teacher, being rational, chooses a day $d\in D$ and reaches her goal,
and this holds even without any announcement. 

\medskip

(B)  $\alpha=\sigma\wedge \langle\cT\in\{\Mo,\We,\Fr\}\rangle$. -- 
The teacher says: ``There will be a surprising test, and it will be on one of the days Monday, Wednesday, and Friday.''
In $\Mc^*$, all worlds $(d,B)$ with $d\in D\cap B$ and $d<\max B$ satisfy this formula. It is easy to see that $D^\alpha=\{\Mo,\We,\Fr\}$,
by just noting that $\none\in B$ is always possible.
The variant $\alpha=\sigma\wedge \langle\cT=\Mo\rangle$ gives us the birthday present ``paradox''~\cite{Gardner:1991,Scriven:1951},
which isn't a paradox, since $D^\alpha=\{\Mo\}$ is the obvious solution set. 
 
\medskip

(C)  $\alpha=\sigma\wedge\nec\langle\cT\in D\rangle$. -- This means that it is guaranteed that only days $d$ in $D$ can be chosen by the teacher as $\rT$.
Technically, only equivalence classes that do not have a {\none}-world are considered. This models the ``paradox of the surprising egg''. 
We have $D^\alpha=\{\Mo,\dots,\Th\}$.

\medskip

(D) $\alpha=\nec\langle\cT=\Mo\rangle\vee\nec\langle\cT=\We\rangle\vee\nec\langle\cT=\Fr\rangle$? -- This artifical example is intended to 
show how the modal logic argumentation works in an unusual situation. Looking at the formula, one may have doubts. 
How should the students ``know'' on Monday morning that it is not the $\nec\langle\cT=\Fr\rangle$-world they are in?
Shouldn't tests on Monday or on Wednesday be surprising, but not a test on Friday?
However, none of the three worlds allowed by $\alpha$ satisfies $\sigma$. What is wrong? 
Well, this formula $\alpha$ is illegal as an announcement in that $\sfdash\alpha\to\sigma$ does not hold. 

\medskip

(E) $\alpha=\sigma\wedge(\nec\langle\cT=\Mo\rangle\vee\nec\langle\cT=\We\rangle\vee\nec\langle\cT=\Fr\rangle)$? -- We try to amend the previous example
by just attaching $\sigma$.  
The teacher wishes to have ``surprise'' and $\nec\langle\cT=\Mo\rangle\vee\nec\langle\cT=\We\rangle\vee\nec\langle\cT=\Fr\rangle$. What happens now?
Well, this formula $\alpha$ is S5-inconsistent. Namely, it is S5-equivalent to 
$(\sigma\wedge\nec\langle\cT=\Mo\rangle) \vee (\sigma\wedge\nec\langle\cT=\We\rangle) \vee (\sigma\wedge\nec\langle\cT=\Fr\rangle)$,
and $\sigma\wedge\nec\langle\cT=d\rangle\equiv_{\text{S5}}\bot$ for all $d\in D$. 
So if the teacher makes this announcement, or indicates that it is her goal to make it happen, the students will recognize it as unsatisfiable.

\begin{remark}\label{rem:inconsistent}
And what should the students do if the teacher makes an announcement $\alpha$ that is S5-inconsistent?
It cannot be used for drawing conclusions, since every formula can be proved from it. 
There is no choice: The students must discard this information completely. 
They cannot even assume the teacher has announced anything.
There could be a test or there could be no test, and if there is a test on some day $d$, it will be surprising.
This is the situation of Example~\ref{examples:A:E}(A), and  
clearly, it is the most uninteresting of all, and it is not the one intended by the story. 
\end{remark}

\smallbox{\paragraph{Acknowledgments.}
I am grateful to quite a few people who listened to or read 
through the ideas that lead to this paper at various stages of maturity,
and pushed it on by pointing out deficiencies, by hints and remarks. 
Uwe R\"osler told me about the paradox. 
Dietrich Kuske and Reinhard Wilhelm were kept busy over a lengthy period of time. 
Thanks are due also to the following: 
Matthias Dietzfelbinger, Bill Gasarch, Oliver Deiser, Timothy Chow, A. Luber, A. and E. Flachs,
Hermann Wilhelm, Stefan Walzer, 
the Research Seminar of the Institute for Theoretical Computer Science 
and the the Seminar for Discrete Mathematics at the Technische Universit\"at Ilmenau.
Anonymous referees pointed out serious deficits in an earlier version.
}

\bibliographystyle{plain} 
\bibliography{refs-SIGACT} % Entries are in the refs.bib file
\end{document}